\documentclass{amsart}

\usepackage{ amsmath,amssymb,amsbsy,amsfonts, latexsym,amsopn,amstext, amsxtra,euscript,amscd, hyperref}
\usepackage{graphicx}
\usepackage{cite}

\usepackage{amsrefs}

\usepackage{url}

\usepackage{xy} \input xy \xyoption{all}

\newcommand\myurl[1]{\url{#1}}

\BibSpec{webpage}{%
  +{}{\PrintAuthors} {author}
  +{,}{ \textit} {title}
  +{}{ \parenthesize} {date}
  +{,}{ \myurl} {myurl}
  +{,}{ } {note}
  +{.}{ } {transition}
}

\newtheorem{thm}{Theorem}

\newtheorem{lem}{Lemma}
\newtheorem{rem}{Remark}
\newtheorem{cor}{Corollary}
\newtheorem{exa}{Example}
\newtheorem{prob}{Problem}
\theoremstyle{defn}
\newtheorem{defn}{Definition}

\def\T3{\mathcal T_3}

\def\cR{\mathcal R}

\def\l{\lambda}
\def\a{\alpha}

\def\iso{\equiv}

\def\Z{\mathbb Z}
\def\Q{\mathbb Q}
\def\C{\mathbb C}

\def\L{\mathcal L}
\def\H{\mathcal H}
\def\M{\mathcal M}
\def\N{\mathcal N}

\def\l{\lambda}

\def\a{\alpha}
\def\p{\mathfrak p}
\def\u{\mathfrak u}

\def\iso{{\, \cong\, }}
\def\<{\langle}
\def\>{\rangle}
\def\X{\mathcal X}
\def\Aut{\mathrm{Aut}}
\def\bAut{\overline {\mathrm{Aut}}}

\def\ss{\mathfrak s}

\def\D{\Delta}
\def\ss{\mathfrak s}
\def\X{\mathcal X}
\def\u{\mathfrak u}
\def\L{\mathcal L}
\def\S{\mathcal S}

\def\ch{\mbox{\rm char } }

\def\J{\tilde J}

\begin{document}

\title{Some remarks on the hyperelliptic moduli of genus 3}

\author{T. Shaska}
\address{Department of Mathematics \& Statistics, Oakland University, Rochester, MI, 48309.}

\email{shaska$@$oakland.edu}

\thanks{The author was supported in part by NSA Grant \#000000.}

%    General info
\subjclass[2000]{Primary 54C40, 14E20; Secondary 46E25, 20C20}

%\date{January 1, 2001 and, in revised form, June 22, 2001.}


\keywords{invariants, binary forms, genus 3, algebraic curves}

\begin{abstract} 
In 1967, Shioda  \cite{Shi1} determined the ring of invariants of binary octavics and their syzygies  using the symbolic method. We discover that the syzygies determined in \cite{Shi1} are incorrect.  In this paper, we compute the correct equations among the invariants of the binary octavics and give necessary and sufficient conditions for two genus 3 hyperelliptic curves to be isomorphic over an algebraically closed  field $k$,   $\ch k \neq 2, 3, 5, 7$.  For the first time, an explicit equation of the hyperelliptic moduli for genus 3 is computed in terms of absolute invariants. 
%Furthermore,  relations among such invariants are determined for all subloci of curves with a fixed automorphism group.  
%Such results will be helpful to mathematicians and cryptographers who work with genus 3 hyperelliptic curves. 
\end{abstract}

\maketitle

\section{Introduction}
%&&&&&&&&&&&&&&&&&&&&&&&&&&&&&&&&&&&&&&&&&&&&&&&&&&&&&&&&&&&&&&&&

Let $k$ be  an algebraically closed field.   A binary  form  of degree $d$ is  a homogeneous  polynomial $f(X,Y)$ of  degree $d$ in two  variables over $k$.  Let $V_d$ be the $k$-vector space of binary  forms of degree $d$.  The group $GL_2(k)$ of  invertible $2  \times 2$  matrices over  $k$  acts on  $V_d$  by coordinate  change. Many  problems  in algebra  involve properties  of binary forms  which are invariant under these  coordinate changes. In particular, any  hyperelliptic genus $g$ curve over  $k$ has a  projective equation of the form $Z^2Y^{2g} = f(X,Y)$,  where $f$ is  a binary form  of degree  $d=2g+2$  and  non-zero discriminant. Two  such curves  are isomorphic  if and  only if  the corresponding  binary forms are conjugate under $GL_2(k)$. Therefore the moduli space $\mathcal H_g$ of hyperelliptic genus $g$ curves  is the affine  variety whose  coordinate ring  is the  ring of $GL_2(k)$-invariants in the coordinate ring  of the set of elements of $V_d$ with non-zero discriminant.  It is well known that the moduli spaces $\H_g$ of hyperelliptic curves of genus $g$, $g\neq 4$, are all rational varieties, i.e. isomorphic to a purely transcendental extension field $k(t_1, \dots , t_r)$; see Igusa \cite{Ig}, Katsylo \cite{Ka}.

Generators for this  and similar invariant rings in  lower degree were constructed by  Clebsch, Bolza  and others in  the last  century using complicated symbolic calculations.  For the case of sextics,  Igusa \cite{Ig} extended  this to
algebraically closed  fields of  any characteristic using  difficult techniques  of algebraic  geometry.  For a modern treatment of the degree six case see \cite{vishi}.

The case of binary octavics has been first studied during the 19th century by  von Gall \cite{vG} and  Alagna \cite{Al, Al1}. Shioda in his thesis \cite{Shi1}  determined the  structure  of the ring of invariants $\cR_8$, which turns out to be generated by nine $SL(2, k)$-invariants $J_2,\cdots, J_{10}$  satisfying five  algebraic relations.  He  computed explicitly these five syzygies, and determined  the corresponding syzygy-sequence and therefore the structure of the ring $\cR_8$; see Shioda \cite{Shi1}.

This paper started as a project to implement an algorithm which determines if two genus 3 hyperelliptic curves are isomorphic over $\C$.  According to Shioda \cite[Thm. 5]{Shi1}; two genus 3 hyperelliptic curves are isomorphic if and only if the corresponding 9-tuples  $(J_2, \dots , J_{10})$ are equivalent, satisfying  five syzygies $R_i (J_2, \dots , J_{10})=0$, for $i=1, \dots , 5$ and non-zero discriminant $\D \neq 0$. While trying to implement the syzygies $R_i (J_2, \dots , J_{10})=0$, for $i=1, \dots , 5$  we discovered that they are not satisfied for a generic octavic.  Hence, such algebraic relations in terms of $J_2, \dots , J_{10}$ are incorrect as stated in \cite{Shi1}; cf Example~\ref{exa_1}.

Indeed, if you take any random binary octavics then its invariants will not satisfy the Shioda's relations.  Since the results in \cite{Shi1} do not hold, then one needs to determine explicitly the algebraic relations between the invariants in order to have an explicit description of the ring of invariants $\cR_8$ and its field of fractions $\S_8$. This will be our goal for the rest of this paper.

In section 2, we give some basic preliminaries on invariants of binary forms. In section 3, we define the main invariants of binary octavics via transvectants. The definitions are the same as used by classical invariant theorists, however, we scale be a constant factor in order to work with primitive polynomials with integer coefficients.  We show  an example of a binary form which does not satisfy the syzygies as claimed in \cite{Shi1}; see Example~\ref{exa_1}. Furthermore, we determine the algebraic relations between the invariants $J_2, \dots , J_{10}$. Such algebraic relations determine the ring of invariants $\cR_8$.

From the basic $SL(2, k)$-invariants $J_2, \dots , J_8$ we define six $GL(2, k)$-invariants 
\[
t_1:= \frac {J_3^2} {J_2^3}, \quad 
t_2:= \frac {J_4} {J_2^2}, \quad 
t_3:= \frac {J_5} {J_2\cdot J_3}, \quad 
t_4:= \frac {J_6} {J_2\cdot  J_4}, \quad 
t_5:= \frac {J_7} {J_2 \cdot J_5}, \quad 
t_6:= \frac {J_8} {J_2^4},
\]
which we call absolute invariants.  There is an algebraic relation \[ T(i_1, \dots , t_6)=0\] that such invariants satisfy, computed for the first time. Shioda in his paper talked about this relation but never attempted to compute it. 
It has total degree     14,    degrees 5, 10, 6, 6, 5, 5 in $t_1, \dots , t_6$ respectively,   and  has 25 464 monomials. 
The field of invariants $\S_8$ of binary octavics is $\S_8= k(t_1, \dots , t_6),$  where $t_1, \dots , t_6$ satisfy the equation $T(t_1, \dots , t_6)=0$. Hence, we have an explicit description of the hyperelliptic moduli $\H_3$.  A birational parametrization of this variety seems out of reach computationally.

%Finally, we have the desired result: two genus 3 hyperelliptic curves are isomorphic if and only if they have the same invariants $t_1, \dots , t_6$.  Moreover, evry point in the variety given by the equation $T (i_1, \dots , t_6)=0 $ correspond to an equivalence class of a genus 3 hyperelliptic curve. 

All of our results are implemented in a Maple package and made available at \cite{homepage}.  Such results will be helpful in the arithmetic of genus 3 hyperelliptic curves.  
The computation of the equation in \eqref{Sh3} makes now possible to describe the subloci of $\H_3$ in terms of the $t_1, \dots , t_6$ invariants and other problems on genus 3 hyperelliptic curves 
as described in  \cite{g_sh, GSS, sh_03, sh_04,  sh_05, issac, ajm_sh1, sh_2000, gen3, sh_thompson, math_comp} among others.

%************************************************************
\section{Preliminaries on invariants of binary forms}
%************************************************************

In this section we define the action of $ GL_2(k)$ on the space of binary forms and discuss the basic notions of their invariants. Most of this section is a summary of section 2 in \cite{vishi}. Throughout this section $k$ denotes an algebraically closed field.

\subsection{Action of $GL_2(k)$ on binary forms.}
%*******************************************************

Let $k[X,Y]$  be the  polynomial ring in  two variables and  let $V_d$ denote  the $(d+1)$-dimensional subspace  of  $k[X,Y]$  consisting  of homogeneous polynomials.
\begin{equation}
\label{eq1} f(X,Y) = a_0 X^d + a_1X^{d-1}Y + ... + a_dY^d
\end{equation}
of  degree $d$. Elements  in $V_d$  are called  \textbf{binary  forms} of degree $d$. We let $GL_2(k)$ act as a group of automorphisms on $ k[X,Y] $ as follows: if
$$ g =
\begin{pmatrix} a &b \\ c & d
\end{pmatrix}
\in GL_2(k) $$ then
\begin{equation}
g \, \begin{pmatrix} X\\ Y \end{pmatrix} =
\begin{pmatrix} aX+bY\\ cX+dY \end{pmatrix}
\end{equation}
This action of $GL_2(k)$  leaves $V_d$ invariant and \textbf{acts irreducibly} on $V_d$.

\begin{rem}
It is well  known that $SL_2(k)$ leaves a bilinear  form (unique up to scalar multiples) on $V_d$ invariant.    This form is symmetric if $d$ is even and skew symmetric if $d$ is odd.
\end{rem}

Let $A_0$, $A_1$,  ... , $A_d$ be coordinate  functions on $V_d$. Then the coordinate  ring of $V_d$ can be   identified with $ k[A_0  , ... , A_d] $. For $I \in k[A_0, ... , A_d]$ and $g \in GL_2(k)$, define $I^g \in k[A_0, ... ,A_d]$ as follows
\begin{equation}
\label{eq_I} {I^g}(f) = I(g(f))
\end{equation}
for all $f \in V_d$. Then  $I^{gh} = (I^{g})^{h}$ and \eqref{eq_I}  defines an action of $GL_2(k)$ on $k[A_0,
... ,A_d]$.
%************************************
\begin{defn}
Let  $\cR_d$  be the  ring of  $SL_2(k)$ invariants  in $k[A_0, ... ,A_d]$, i.e., the ring of all $I \in
k[A_0, ... ,A_d]$ with $I^g = I$ for all $g \in SL_2(k)$.
\end{defn}

Note that if $I$ is an invariant, so are all its homogeneous components. So $\cR_d$ is  graded by the  usual  degree  function on  $k[A_0, \dots ,A_d]$.

Since $k$ is algebraically closed, the binary form $f(X,Y)$ in Eq.~\eqref{eq1} can be factored as
 \begin{equation}  f(X,Y)  = (y_1  X  -  x_1 Y)...  (y_d  X  - x_d  Y)
=  \prod_{1 \leq  i \leq  d} det
\begin{pmatrix}
X&x_{i}\\Y&y_i\\
\end{pmatrix}
\end{equation}
The points  with homogeneous coordinates $(x_i, y_i)  \in \mathbb P^1$ are  called the \textbf{ roots  of the  binary  form} in Eq.~ \eqref{eq1}.  Thus  for $g  \in GL_2(k)$ we have

\begin{center}
$g(f(X,Y))  = (det(g))^ d  (y_1^{'}  X -  x_1^{'}  Y)...(y_d^{'} X  - x_d^{'} Y)$.
\end{center}
where
\begin{equation}
 \begin{pmatrix}   x_i^{'}  \\  y_i^{'}  \end{pmatrix}
= g^{-1}
\begin{pmatrix} x_i\\ y_i \end{pmatrix}
\end{equation}

%*****************************************
%\subsection{The Null Cone of $V_d$}
%********************************************************

\begin{defn}
The  \textbf{nullcone}  $\N_d$ of  $V_d$  is the  zero set  of all  homogeneous elements in $\cR_d$ of positive    degree
\end{defn}

\noindent The notion of {\it nullcone} was first used by Hilbert; see \cite{H}. Next we define the {\it   Reynold's operator} on $k[A_0, \dots , A_d]$.

\begin{lem}
\label{lem1} Let $char(k)  = 0$  and $\Omega_s$  be the subspace  of $k[A_0,  ... , A_d]$ consisting of  homogeneous elements  of degree $s$. Then there is a $k$-linear  map
$$R :  k[A_0, ... ,  A_d] \to \cR_d$$
with the following properties:

(a) $R(\Omega_s) \subseteq  \Omega_s$ for all $s$

(b) $R(I) = I$ for all $I \in \cR_d$

(c) $R(g(f)) = R(f)$ for all $f \in k[A_0, ... , A_d]$
\end{lem}

\proof $\Omega_s$ is a  polynomial module of degree $s$  for $SL_2(k)$. Since $SL_2(k)$  is linearly  reductive in $char(k)  = 0$,  there exists  a $SL_2(k)$-invariant  subspace  $\Lambda_s$  of $\Omega_s$  such
that $\Omega_s = (\Omega_s \cap  \cR_d) \bigoplus \Lambda_s$. Define
$$R : k[A_0,  ... , A_d] \to  \cR_d$$
such that $R(\Lambda_s)  = 0$ and $R_{|\Omega_s \cap \cR_d} = id$. Then $R$ is $k$-linear and the rest of the    proof is clear from the definition of $R$.

\qed

\noindent The map $R$ is called the \textbf{Reynold's operator}.

\begin{lem}
\label{lem2} Suppose $char(k) = 0$.  Then every maximal ideal in $ \cR_d$ is contained in a maximal ideal of
$k[A_0, ... , A_d]$.
\end{lem}

\proof If $\mathcal I$ is a maximal ideal in $ \cR_d $ which generates the unit ideal of $k[A_0, ... , A_d]$,  then there exist $m_1$, $ m_2$, ... ,$m_t \in \mathcal I$ and $f_1$, $f_2$, ... ,$f_t \in k[A_0, ... , A_d]$
such that
\begin{center}
$1 = m_1 f_1 + ... +m_t f_t$
\end{center}
Applying the Reynold's operator to the above equation we get
\begin{center}
$1 = m_1 R(f_1) + ... + m_t R(f_t)$
\end{center}
But  $R(f_i)  \in \mathcal  R_d$  for all  $i$.  This  implies $1  \in \mathcal I$, a contradiction. \qed

\begin{thm}\label{thm1} [Hilbert's Finiteness  Theorem]  Suppose $char(k)  = 0$.  Then $\cR_d$ is finitely  generated over $k$. \label{Hilbert1}
\end{thm}

\proof Let $\mathcal I_0$ be the ideal in $k[A_0, ... ,A_d]$ generated by all homogeneous invariants of  positive degree. Because $k[A_0, ... ,A_d]$ is Noetherian,  there exist finitely many  homogeneous elements
$J_1$, ... $J_r$  in  $\mathcal  R_d$   such  that  $\mathcal  I_0  =  (J_1, ... ,J_r)$. We prove $\cR_d =  k[J_1, ... , J_r]$.  Let $J \in \cR_d$  be homogeneous  of degree $d$.  We prove $J  \in k[J_1, ... , J_r]$
using induction on $d$.  If $d = 0$, then $J \in k \subset k[J_1, ... , J_r]$. If $d > 0$, then
\begin{equation}
\label{eq6} J = f_1 J_1 + ... + f_r J_r
\end{equation}
with $f_i  \in k[A_0, ... ,  A_d]$ homogeneous and $deg(f_i)  < d$ for all $i$. Applying the Reynold's  operator to Eq.~\eqref{eq6} we have
$$J = R(f_1) J_1 + ... + R(f_r) J_r$$
then by  lemma \ref{lem1}, we have  $R(f_i)$ is  a homogeneous element  in $\cR_d$ with $deg(R(f_i))  < d $    for all $i$  and hence by induction  we have $R(f_i) \in k[J_1, ...  , J_r]$ for all $i$. Thus $J  \in k[J_1, \dots  , J_r]$. \qed

If $k$ is of arbitrary characteristic, then $SL_2(k)$ is geometrically reductive,  which is a weakening of  linear reductivity; see Haboush \cite{Ha}. It suffices to prove Hilbert's finiteness theorem in any  characteristic; see Nagata \cite{Na}. The following theorem is also due to Hilbert \cite{H}.

%******************************************

\begin{thm}\label{thm2} Let $I_1$,  $I_2$, \dots ,  $I_s$ be homogeneous  elements in  $ \mathcal R_d$ whose  common zero set equals the null cone $\mathcal N_d$. Then $ \cR_d$  is finitely  generated as a module over  $k[I_1, \dots , I_s]$.
\end{thm}

\proof (i) $char(k)  = 0$: By theorem \ref{thm1},   we have $\cR_d  = k[J_1, J_2, \dots  ,  J_r]$ for some
homogeneous invariants  $J_1, \dots  , J_r$. Let $\mathcal I_0$  be the maximal ideal  in $ \mathcal  R_d$
generated by all homogeneous elements  in $ \cR_d$ of  positive degree. Then the theorem follows if $I_1,
\dots , I_s$ generate an ideal $\mathcal I$ in $ \cR_d$ with $rad (\mathcal I) = \mathcal I_0$. For if this
is the case, we have an integer $q$ such that
\begin{equation}
\label{eq7}
 J_i ^{q} \in \mathcal I \hspace{0.08in} \forall i
\end{equation}
Set $S:= \{ J_1^{i_1} J_2^{i_2} \dots J_r^{i_r} : 0 \leq i_1, \dots , i_r < q \}$. Let $\M$ be  the $k[I_1,
\dots ,  I_s]$-submodule in $\cR_d$  generated by $S$.  We prove $\cR_d  = \M$. Let $J  \in \mathcal  R_d$ be
homogeneous.  Then $J  = J^{'}  + J^{''}$  where  $J^{'} \in  \M$,  $J^{''}$  is a  $k$-linear combination of
$J_1^{i_1} J_2 ^{i_2}  ... J_r ^{i_r}$ with  at least one $i_{\nu} \geq  q$ and $deg(J) = deg(J^{'})  =
deg(J^{''})$. Hence \eqref{eq7} implies $J^{''} \in \mathcal I$ and so we have
\[J^{''} = f_1 I_1 + \dots  + f_s I_s,\]
where $f_i \in \cR_d$ for all $i$. Then $deg(f_i) < deg(J^{''}) = deg(J)$ for all $i$. Now by induction on  degree of $J$ we may assume $f_i \in \mathcal M$ for all $i$. This implies $J^{''} \in \mathcal M$ and  hence
$J \in  \mathcal  M$.  Therefore  $\mathcal M  =  \mathcal R_d$.  So  it only  remains  to prove   $rad(\mathcal  I)  =  \mathcal I_0$. This  follows from  Hilbert's Nullstellensatz and  the following claim.

\smallskip

\noindent {\bf Claim:}  $\mathcal I_0$ is the only maximal ideal containing $I_1$, \dots  , $I_s$.

Suppose $\mathcal  I_1 $ is  a maximal ideal  in $ \mathcal  R_d$ with $I_1, ...  ,I_s \in \mathcal I_1$.  Then from Lemma  2 we  know there exists a  maximal ideal  $\mathcal J$  of $k[A_0, ...  , A_d]$  with $\mathcal I_1 \subset \mathcal J$.  The point in $V_d$ corresponding to $\mathcal  J$ lies  on the  null  cone $\mathcal  N_d$ because  $I_1$, ...  ,$I_s$ vanish  on  this point.  Therefore  $\mathcal I_0  \subset \mathcal J$,  by definition of  $\mathcal N_d$. Therefore $\mathcal J \cap \cR_d$ contains both the maximal
ideals $\mathcal I_1$ and $\mathcal I_0$. Hence $\mathcal I_1 = \mathcal J  \cap \cR_d = \mathcal I_0$.

\smallskip

(ii)  $char(k) = p$: The same proof  works if lemma \ref{lem2} holds. Geometrically  this means the morphism
$\pi : V_d  \to V_d$ // $SL_2(k)$ corresponding to the  inclusion $\cR_d \subset k[A_0, ... , A_d]$ is surjective. Here $ V_d$ // $SL_2(k)$ denotes the affine variety corresponding  to the  ring $\cR_d$  and is
called the \textbf{categorical  quotient}. $\pi$ is surjective  because $SL_2(k)$ is geometrically  reductive.  The proof  is by  reduction modulo  $p$;  see Geyer \cite{Ge}. 

\qed

%****************
\subsection{Hyperelliptic curves of genus 3}
%*********************************************
In this section we want to use the projective equivalence of binary octavics in order two give conditions that two hyperelliptic curves of genus 3 are isomorphic.  

Denote a binary form of order $2g+2$ by 
\[ f (X, Y) = \sum_{i=0}^{2g+2} a_i X^i\, Y^{2g+2-i}\]
To each $f(X, Y)$ with no multiple roots we associate the non-singular hyperelliptic curve $C_f$ with affine equation $Z^2= f(X, 1)$.  Every hyperelliptic curve of genus $g$ is obtained this way. 

Two hyperelliptic curves $C_f$ and $C_h$ are birationally equivalent  if and only if $f(X, Y)$ and $h(X, Y)$ are projectively equivalent, i.e., there exists a $\tau \in SL_2 (k)$ and $\l \in k\setminus \{0\}$ such that $f^\tau = \l \cdot h$.  

Let $\D_f$ denote the discriminant of the polynomial $f(X, 1)$.  It is an invariant of degree $2 (2g+1)$.  When $g=3$ then the discriminant has degree 14 and is given as a polynomial in $J_2, \dots , J_8$.

%****************
\section{Projective invariance of binary octavics.}

Throughout this section $\ch (k) \neq 2, 3, 5, 7$.

%*******************************************
\subsection{Covariants and invariants of binary octavics}
%*************************************************
We will use the symbolic method of classical theory to construct covariants of binary octavics. They were
first constructed by van Gall who showed that there are 70 such covariants; see von Gall \cite{vG}. First we  recall
some facts about the symbolic notation. Let
\[f(X,Y):=\sum_{i=0}^n
\begin{pmatrix} n \\ i
\end{pmatrix}
a_i X^{n-i} \, Y^i, \quad  and \quad g(X,Y) :=\sum_{i=0}^m
  \begin{pmatrix} m \\ i
\end{pmatrix}
b_i X^{n-i} \, Y^i
\]
be binary forms of  degree $n$ and $m$ respectively. We define the $r$-{\it transvection}
\[(f,g)^r:= \frac {(m-r)! \, (n-r)!} {n! \, m!} \, \,
\sum_{k=0}^r (-1)^k
\begin{pmatrix} r \\ k
\end{pmatrix} \cdot
\frac {\partial^r f} {\partial X^{r-k} \, \,  \partial Y^k} \cdot \frac {\partial^r g} {\partial X^k  \, \,
\partial Y^{r-k} },
\]
see Grace and Young \cite{GY} for  details.

The following result gives relations among the invariants of binary forms and it is known as the Gordon's formula.  It is the basis for most of the classical papers on invariant theory.  

\begin{thm}[Gordon]
Let $\phi_i$, $i=0, 1, 2$  be covariants of order $m_i$ and $e_i$ be three non-negative integers such that $e_i+e_j \leq m_k$ for distinct $i, j, k$.  The following is true:
%
%\begin{Small}
\[  
\sum_i  \frac  {C_i^{e_1} \cdot  C_i^{m_1-e_0-e_2 }}  {C_i^{m_0+m_1+1-2e_2-i}   } \, \left(  \left(\phi_0 \, \phi_1 \right)^{e_2+1}, \phi_2  \right)^{e_0+e_1-i}  
= \sum_i  \frac  {C_i^{e_2} \cdot  C_i^{m_2-e_0-e_1 }}  {C_i^{m_0+m_2+1-2e_1-i}   } \, \left( \left(\phi_0 \, \phi_2 \right)^{e_1+1}, \phi_1      \right)^{e_0+e_2-i}, 
\]
%\end{Small}
where $e_0=0$ or $e_1+e_2=m_0$.  
\end{thm}

This result has been used by many XIX century mathematicians to compute algebraic relations among invariants, most notably by Bolza for binary sextics and by Alagna for binary octavics.  It provides algebraic relations among the invariants in a very similar manner that the Frobenious identities do for theta functions of hyperelliptic curves. Whether there exists some explicit relation among both formulas seems to be unknown.  

For the rest of this paper $f(X,Y)$ denotes a binary octavic as below:
\begin{equation}
f(X,Y) =   \sum_{i=0}^8 a_i X^i Y^{8-i} = \sum_{i=0}^8
\begin{pmatrix} n \\ i
\end{pmatrix}    b_i X^i Y^{n-i}
\end{equation}
where $b_i=\frac {(n-i)! \, \, i!} {n!} \cdot a_i$,  for $i=0, \dots , 8$. We define the following
covariants:
\begin{equation}
\begin{split}
&g=(f,f)^4, \quad k=(f, f )^6, \quad h=(k,k)^2, \quad m=(f,k)^4, \quad  n=(f,h)^4, \quad p=(g,k)^4, \quad q=(g, h)^4.\\
\end{split}
\end{equation}

\noindent Then, the following 
\begin{equation}\label{def_J}
\begin{aligned}
&     J_2= 2^2 \cdot 5 \cdot 7 \cdot (f,f)^8,                    &   \qquad   &  J_3= \frac 1 3 \cdot  2^4 \cdot 5^2 \cdot 7^3 \cdot (f,g )^8, \\
&  J_4= 2^9 \cdot 3 \cdot 7^4 \cdot (k,k)^4,                      &    \qquad      &   J_5= 2^9 \cdot 5 \cdot 7^5 \cdot (m,k)^4,  \\
&   J_6 = 2^{14} \cdot 3^2 \cdot 7^6 \cdot (k,h )^4,              &    \qquad      &      J_7= 2^{14} \cdot 3 \cdot 5 \cdot 7^7 \cdot (m,h )^4,  \\
& J_8= 2^{17} \cdot 3 \cdot 5^2 \cdot 7^9 \cdot   (p,h)^4,        &    \qquad       &   J_9= 2^{19} \cdot 3^2 \cdot 5 \cdot 7^9 \cdot   (n,h)^4, \\
&  J_{10}=   2^{22} \cdot 3^2 \cdot 5^2 \cdot 7^{11} (q,h)^4      &  \qquad    &       \\
\end{aligned}
\end{equation}
are $SL_2(k)$- invariants.     Notice that we are scaling such invariants up to multiplication by a constant for computational purposes only.  We display only the first two of such invariants to avoid any confusion in the definitions

\[
\begin{split}
J_2  = & 280\,a_8a_0-35\,a_7a_1+10\,a_6a_2-5\,a_5a_3+2\,{a_4}^2 \\
J_3  = & 1050\,a_8{a_2}^2+1050\,{a_6}^2a_0+75\,a_6{a_3}^2+75\,{a_5}^2a_2  +12\,{a_4}^3+3920\,a_8a_4a_0-2450\,a_8a_3a_1+735\,a_7a_4a_1\\
& -2450\,a_7a_5a_0 -175\,a_7a_3a_2-110\,a_6a_4a_2-175\,a_6a_5a_1-45\,a_5a_4a_3 \\
\end{split}
\]
 In other words, we take the numerator of the corresponding transvectants since we prefer to work over $\Z$ instead of $\Q$ and then take the primitive part of each invariant. Hence,  we have $J_i \in \Z [a_0, \dots , a_8]$, for $i= 2, \dots , 8$ and $J_i$'s are primitive polynomials. 
 In \cite{Shi1} such scaling is not done and this invariants are homogenous polynomials with coefficients in $\Q [a_0, \dots , a_8]$ and not  primitive. 

\begin{lem}\label{lem_3} For each binary octavic $f(X, Y)$, its invariants defined in Eq.\eqref{def_J}  are primitive homogeneous polynomials $J_i\in \Z [a_0, \dots , a_8]$  of degree $i$, for $i=2, \dots , 10$.
Let $f^\prime=g (f)$, where
\[ g =
\begin{pmatrix} a &b \\ c & d
\end{pmatrix}
\in GL_2(k),  
\]
and denote the corresponding $J_2, \dots , J_{10}$ of $f^\prime$ by $J_2^\prime, \dots , J_{10}^\prime$. Then,
$$J_i^\prime= ( \Delta^4)^i \, J_i$$ where $\Delta=ad-bc$ and $i=2, \dots , 10$.
\end{lem}

\proof The first claim is immediate from the definition of the covariants and invariants.  Let $f$ and
$f^\prime$ be two binary octavics as in the hypothesis. One can check the result computationally. \qed

\begin{rem}
There are 68 invariants defined this way as discovered by van Gall \cite{vG, vG1} in 1880.  Indeed, van Gall claimed 70 such invariants, but as discovered in XX-century there are only 68 of them.  
Perhaps, one that needs to be mentioned is $J_{14}$ which is the discriminant of the binary octavic. 

In a couple of papers in 1892 and 1896 R. Alagna determined the algebraic relations among such invariants; see \cite{Al, Al1} for details.  All these works have computational mistakes and are almost impossible to check. 
\end{rem}

Next we want to show that the ring of invariants $\cR_8$ is finitely generated as a module over $k[J_2,
\dots,J_7]$. First we need some auxiliary lemmas.

\begin{lem} If $J_i=0$, for $i=2, \dots 7$, then the $f(X,Y)$ has a  multiple root.
\end{lem}

\proof Compute $J_i=0$, for $i=2, \dots 7$. These equations imply that
\[Res( f(X,1), f^\prime (X, 1), X)=0,\]
where $f^\prime $ is the derivative of $f$. This proves the lemma. \qed

% ************************
\begin{thm} \label{thm_5}   The following hold true for any octavic.

i)   An octavic has a root of multiplicity exactly four if and only if the basic invariants take the form
\begin{equation}\label{J_i}
\begin{split}
J_2 &= 2  \cdot r^2, 
\quad       J_3=  2^2 \cdot 3 \cdot  r^3,  
\quad J_4 = 2^6 \cdot   r^4, 
\quad J_5 = 2^6   \cdot r^5,\\
 &J_6=  2^9 \cdot    r^6, 
 \quad J_7= 2^9 \cdot   r^7, 
 \quad J_8=  2^{11}\cdot 3^2 \cdot        r^8,
\end{split}
\end{equation}
for some  $ r \neq 0 $.  Moreover, if the octavic has equation 
\[ f(x, y)=x^4 (a x^4+b x^3 y+c x^2 y^2+d x y^3+e y^4), \]
then  $r=e$.

ii) \label{lem_root} An octavic has a root of multiplicity 5 if and only if
\[ J_i=0, \ \  for \ \ i=2, \dots , 8.\]

\end{thm}

\proof i)   Let
$$f(X,Y) = a_0 X^8 + a_7 X^7 Y + \dots + a_8 Y^8$$
be an octavic with a root of multiplicity four. Let this root be at $(1,0)$. Then,
\[f(X,Y) = (a_4 X^4 + a_3 X^3 Y + a_2 X^2 Y^2+ a_1 X Y^3 +a_0 Y^4) X^4\]
 Thus,  for $r=a_4$,  $\, J_{i}$ for $ i = 2, \dots , 8$ are as claimed.

Conversely assume that Eq.~\eqref{J_i} holds. Then, we have a multiple root. We assume the multiple root is at
$(1,0) $. If this is the only root then $r=0$. Thus, there is at least one more root.  We assume the  other
root is $ ( 0,1) $. Then the octavic takes the form
\begin{equation}
f(X,Y)=a_2 X^6 Y^2 + a_3 X^5 Y^3 + a_4 X^4 Y^4 + a_5 X^3 Y^5 + a_6 X^2 Y^6 + a_7 X Y^7
\end{equation}
and \eqref{J_i} becomes a system of six equations. We eliminate $a_2, a_3$ to get that $a_5=0 $ or $a_4=r$.
If $a_4=r$ and $a_5\neq 0$ then $a_2=a_3=0$ and $(1,0)$ is a root of multiplicity four. If $a_5=0$ then from
the system we get $a_2=0$ or $a_6=0$. In both cases we have a root of multiplicity four. \qed

ii)  Suppose $(1, 0)$ is a root of multiplicity 5. Then, as in  previous lemma  we can take $a_8= a_7 = a_6=a_5=a_4= 0 $. Then by a lemma of Hilbert \cite{H} or  by simple computation we have
 $J_i= 0$, for $i=2, \dots , 7$.

For the converse, since $J_{14}= 0$, there is a multiple root. If there is no root other than the multiple
root, we are done. Otherwise, let the multiple root be at (1,0) and the other root be at (0, 1). Then as in
the previous lemma, the octavic becomes
\begin{equation}
f(X,Y)=a_2 X^6 Y^2 + a_3 X^5 Y^3 + a_4 X^4 Y^4 + a_5 X^3 Y^5 + a_6 X^2 Y^6 + a_7 X Y^7
\end{equation}

\noindent Compute all $J_2, \dots J_7$. From the corresponding system of equations we can eliminate $a_2, a_3,
a_7$. We have a few cases:
%
%\begin{Small}
\[a_4  \left( -2\,{\it a_4}\,{\it a_6}+{{\it a_5}}^2 \right)
 \left( -34\,{\it a_4}\,{\it a_6}+15\,{{\it a_5}}^2 \right)
 \left( 5476\,{{\it a_6}}^2{{\it a_4}}^2+2025\,{{\it a_5}}^4-
6780\,{\it a_4}\,{{\it a_5}}^2{\it a_6} \right) =0\]
%\end{Small}
%
Careful analysis of each case leads to the existence of a root of multiplicity 5. The proof is computational
and we skip the details. \qed

\begin{rem}
An alternative proof of the above can provided using the $k$-th subresultants of $f$ and its derivatives. Two
forms have $k$ roots in common if and only if the first  $k$ subresultants vanish.
%, for details see \cite{Boc}. 
This is equivalent to $J_2= \dots = J_7=0$.
\end{rem}

\subsection{The Null Cone of $V_8$ and Algebraic Dependencies}
%***************************************************************

\begin{thm}\label{thm_6}
$\cR_8$ is finitely generated as a module over $k[J_2, \dots,J_7]$.
\end{thm}

\proof By Theorem 2 we only have to prove $\mathcal N_8 = V(J_2, \dots ,J_7)$. For $\l \in k^{*}$, set
$$ g(\l):=
\begin{pmatrix}
\l ^{-1} &0\\ 0 &\l
\end{pmatrix}
$$
Suppose $J_2, \dots , J_7$ vanish on an octavic $f\in V_8$. Then  we know from Theorem~\ref{lem_root} that $f$
has a root of multiplicity at least 5. Let this multiple root be $(1,0)$. Then $f$ is of the form
$$f(X,Y) = (a_5 X^3 + a_6 X^2 Y + a_7 X Y^2 + a_8 Y^3)\,  Y^5$$
If  $I \in  \cR_8$ is  homogeneous  of degree  $s>0$,  then
$$I(f^{g(\l)})=\l^{2s} I(a_5 X^3 Y^5 +\l^2 a_6 X^2 Y^6
+\l^3 a_7 X Y^7 +\l^4 a_8 Y^8 )$$
Thus $I(f^{g(\l)})$ is a polynomial in   $\l$  with   no  constant   term.  But   since  $I$   is an
$SL_2(k)$-invariant,  we  have  $I(f^{g(\l)})=I(f)$  for  all $\l$. Thus $I(f)=0$. Then, ${\mathcal N_8}=V
(J_2, J_3, J_4, J_5, J_6, J_7)$. This completes the proof.

\qed

The above lemma is proven by Shioda in a more computational way using the symbolic method; see below for more details.

\begin{cor}
$J_2, \dots , J_7$ are algebraically independent over $k$ because $\cR_8$ is the coordinate ring of the     5-dimensional variety $V_8$//$SL_2(k)$.
\end{cor}

\subsubsection{Shioda's computations}
The algebraic relations between $J_2, \dots , J_{10}$ were computed by Shioda in \cite{Shi1} using the symbolic method.  However, we could not confirm the correctness of such results with our computations.  For a binary octavic 
\[ f (X, Y) = \sum_{i=0}^8 a_i X^i Y^{8-i}, \] 
Shioda invariants are defined as 
\[
\begin{split}
\J_2 = &   2\,a_8a_0-16\,a_7a_1+56\,a_6a_2-112\,a_5a_3+70\,{a_4}^2\\
\J_3 = &  {\frac 9{392}}\,a_8{a_2}^2+{\frac 9{392}}\,{a_6}^2
a_0+{\frac 9{5488}}\,a_6{a_3}^2+{\frac 9{5488}}\,{a_
5}^2a_2+{\frac 9{560}}\,a_7a_4a_1-{\frac 3{
56}}\,a_7a_5a_0 \\
& -{\frac 3{784}}\,a_7a_3a_2-{ \frac {33}{13720}}\,a_6a_4a_2-{\frac 3{784}}\,a_6a_{5}a_1-{\frac {27}{27440}}\,a_5a_4a_3+{\frac 3{35}} \,a_8a_4a_0-  {\frac 3{56}}\,a_8a_3a_1+{\frac 9{34300}}\,{a_4}^3 \\
\end{split} 
\]
Notice that that definition of $J_2$ looks different from that of Shioda \cite[page 1037]{Shi1}, but that is because there $J_2$ is evaluated for $ f (X, Y) =   \sum_{i=0}^8 \, {8 \choose i } \,  a_i X^i Y^{8-i}$.

Now we are ready to show that the syzygies in \cite[Th. 5]{Shi1} are not correct.  Below is an example of a genus 3 hyperelliptic curve with invariants which do not  satisfy Shioda relations.

\begin{exa}\label{exa_1}
Let a genus 3 hyperelliptic curve be given by the equation 
\[ y^2= x^8+x^7+x^6+ x^5+x^4+x^3+x^2+x+1\]
Then, its invariants are 
\[ 
\begin{split}
& \J_2= \frac 9 5, \, \J_3= \frac {81} {2450}, \, \J_4= \frac {837} {1568}, \, \J_5=\frac {2187} {109760}, \, \J_6= - \frac {6885} {43904}, \\
&  \J_7= - \frac {3645 } {1229312 },  \,  \J_8= - \frac {410427} {17210368}, \,   \J_9= \frac {234009}  {172103680 }, \,  \J_{10}= \frac  {5972697}  {860518400 }
\end{split}
\]
Then evaluating all expressions as in Shioda's paper  we have 
\begin{Small}
\[
\begin{aligned}
 A_6 & = =  - \frac {3645} { 9604},                       &         A_7 & =  \frac   {130491}  {439040},                           &    A_8 & =  -  \frac {15261615}  {87808},\\
  B_7 & =  \frac {130491}  {351232},                      &        B_9 & =  \frac {1414989}  {172103680},                          &    B_8 & =  \frac  {143437311}  {21512960}, \\
  C_9 & =  \frac {809753208633}  {376476800000},          &       C_{10} & =  -  \frac {51828148570131}  {150590720000},           & & \\
  D_{10} & =  -   \frac  {19194738471171} {385512243200}, &  A_{16} & =  -  \frac {1097050897751848407621}  {925614895923200000}.  &    \\
\end{aligned}
\]
\end{Small}
Substituting all these values in the first equation of   \cite[Thm. 5]{Shi1} we get the value 
\[ -   \frac  {546607935510034107123}  {462807447961600000} \neq 0. \]
This implies that the relations determined by Shioda are not correct.
\end{exa}

Indeed, if you take any random binary octavics then its invariants will not satisfy the Shioda's relations.  Since the results in \cite{Shi1} do not hold, then one needs to determine explicitly the algebraic relations between the invariants in order to have an explicit description of the ring of invariants $\cR_8$ and its field of fractions $\S_8$. This will be our goal for the rest of this paper.

%************************************************************** 
\subsubsection{Algebraic dependencies among the invariants}
In this section we will determine algebraic relations among the invariants $J_2, \dots , J_8$.  We will use computational algebra techniques such as elimination by resultants, Groebner bases, etc.  Any computer algebra package can be used to reproduce our results.  Once obtained, these results can be easily verified.  All our results are organized in a Maple package and will be freely made available.

% 
%\begin{thm}   The  invariants  $J_2, \dots , J_{10}$ of a general  binary octavic satisfy  the following equations
%
%\begin{equation}\label{system}
%
%\left\{
%\begin{split}
% f_1(J_2, \dots J_{10})   =0 \\
%  f_1(J_2, \dots J_{10})   =0 \\
%   f_1(J_2, \dots J_{10})   =0 \\
%    f_1(J_2, \dots J_{10})   =0 \\
%     f_1(J_2, \dots J_{10})   =0 \\
%\end{split}
%\right.
%\end{equation}
%\end{thm}

%\proof

%The proof of the theorem can follow simply by computing the invariants of a generic binary octavic and substituting them into the equation.  They will satisfy the equation.  

%We briefly describe how we computed these equations.      

Without loss of generality we can assume that the generic binary octavic is given by
\begin{Small}
\begin{equation}\label{ini_form}
\begin{split}
f(X, 1)  &= X(X-1)(X^5 - s_1 X^4 + s_2 X^3 - s_3 X^2 + s_4 X - s_5)\\
&=  X^7-(s_1+1) X^6+(s_2+s_1) X^5-(s_3+s_2) X^4+(s_4+s_3) X^3-(s_4+s_5) X^2+s_5 X\\
\end{split}
\end{equation}
\end{Small}
Denote by
\[ a:=s_1+s_2, \  b:=s_2+s_3, \ c:=s_3+s_4, \ d:=s_4+s_5, \  s:=s_5.\]
Then we have
\[ f(X, 1)=X^7+(-1+b-s+d-c-a)X^6+a X^5-b X^4+c X^3-d X^2+s X\]
We first compute the $J_2, \dots , J_{10}$ for $f(X, 1)$.  
\[
\begin{split}
J_2 (f)=&  -35\,s+10\,d-10\,db+10\,ds-10\, d^2+10\,dc+10\,da-5\,ac+2\, b^2 \\
J_3 (f) =&  -75\,c^2+75\,c^2b-75\,c^2s+75\,c^2d-75\,c^3-75\,c^2a-75\,d\, a^2-12\, b^3+110\,db  -110\,d b^2+110\,dbs-110\,{d }^2b\\
& +110\,dbc  +110\,dba+175\,as-175\,asb+175\,a{s}^2 -175\,asd+175\, asc+175\,\, a^2s-735\,bs+175\,dc+45\,cba, \\
\end{split}
\]
$J_4, \dots , J_8$ are largeer expressions and we do not display them.  

%Naively, one would try to eliminate $a, b, c, d$ and $s$ from the system of equations and obtain algebraic relations in terms of the invariants.  This will not be easy to do with a standard  way via Gr\"obener bases.  
%We painfully use resultants to eliminate each one of the parameters.  In every step the resultant polynomial is factored and the content is cancelled.  
%
%\endproof

Our goal is to express $J_8, J_9, J_{10}$ in terms of $j_2, \dots , J_7$.  Indeed, from Thm.~\ref{thm_6} it is enough to express $J_8$ in terms of $J_2, \dots , J_7$.  Since in \cite{Shi1} the syzygies include expressing $J_9$ and $J_{10}$ in terms of $J_2, \dots , J_7$ we will comment on how that can be done also.  

We have the following system of equations
\begin{equation}\label{system}
\left\{
\begin{split}
F_2:= J_2-J_2(a,b,c,d,s)   =0 \\
F_3:= J_3-J_3(a,b,c,d,s)   =0 \\
F_4:= J_4-J_4(a,b,c,d,s)   =0 \\
F_5:= J_5-J_5(a,b,c,d,s)   =0 \\
F_6:= J_6-J_6(a,b,c,d,s)   =0 \\
F_7:=  J_7-J_7(a,b,c,d,s)   =0 \\
F_8:=   J_8-J_8(a,b,c,d,s)   =0 \\
\end{split}
\right.
\end{equation}

We compute the equation of $J_8$ in terms of $J_2, \dots , J_7$ using the following technique.  Take the resultant with respect to $a$ of the polynomials  $F_i, F_8$, for $i=2, \dots 7$.  Let $G_i:=Res(F_i, F_8, a)$, for $i=2, \dots 7$.  For each resultant we want to factor the result and take the primitive part.  It is exactly this part that is important and it is not usually done by implementations of Grobener basis algorithms.  In many cases the resultant will be factored to a power or will have factors which imply that $J_{14}=0$.  Since we are computing in an integral domain, we  cancel such factors. 

We continue now with the system $G_i:=Res(F_i, F_8, a)$, for $i=2, \dots 7$ and compute the resultants $H_i:=Res (G_i, G_7, b)$ for $i=2, \dots 6$.  Hence, we are left 5 equations and transcendentals $c, d, s$.  Continuing this process we get a degree 8 equation of $J_8$ in terms of the other $J_2, \dots , J_7$,   as expected by Shioda; see \cite[pg. 1044]{Shi1}. 
Its leading monomial has coefficient  $2^2 \cdot 3^{20}\cdot 5^{12}$. Since we are assuming that the characteristic of the field is $ \neq 2, 3, 5, 7$ then we can divide by this coefficient.  Hence, denote the minimal quintic   by  
\[J_8^5 + c_4 J_8^4 + \cdots + c_1 J_8 + c_0 =0.\] 

Since this equation is a homogenous equation of degree 40 in $J_2, \dots , J_8$, then all other coefficients of $J_8$ are homogenous polynomials in $J_2, \dots , J_7$ of degree  8, 16, 24, 32, 40 respectively.  We denote the primitive part of each of these coefficients by  $I_8, I_{16}, I_{24}, I_{32}, I_{40}$.      

For now on we use the following notation  
\[ J_2:=a, \quad  J_3:=b, \quad  J_4:=c, \quad J_5:= d, \quad J_6:=e, \quad J_7:=f, \]
to display the expressions of $I_8, I_{16}, I_{24}, I_{32}, I_{40}$. 
\begin{Small}
\begin{equation*}
\begin{split}
I_8 & = -2^7  7^5 a^4+ 2^2  5^3   7^3  3 a^2 c+ 2^6   3^3   7^2 a b^2 +2^3   5^4  7^2 a e- 2^2  3^5  5^2  7 b d- 3^3  5^4  17 c^2 \\
I_{16} & = 2^2  3^7  5^5  7^5  a^3  d^2+   2^2  3^8  5^6  7\cdot 11  b c^2 d +2^3  5^8  7^4  a^2  e^2  -2\cdot 3^7  5^6  7^3  a  b  c  f+2^9  3^6  7^4 a^2  b^4  -3\cdot 2^3  5^5  7^6\cdot 31  a^3  c  e  -2^2  3^8  5^6  7^3  a  c  d^2\\
&  + 2^3  3^4  5^5  7^2  11^2  b^2  c  e - 2\cdot 3^3  5^7  7^3  13  a^2  c^3- 5717  3^5  5^4  7^2  a  b^2  c^2+2^6  3^6  5^3  7^5  a^3  b  f+2^2  3^8  5^5  7^4  a^2  d  f  +2^7  3^5  5^2  7^6  a^4  b  d\\
&  -3^8  5^7  7^3  a  f^2 -2^3  3^{10}  5^4  7^2  b^2  d^2-2^5  3^7  5^2  7^2 19  b^4  c-2^2  3^7  5^7  7^2  d^2  e+2^{11}  7^{10}  a^8+  2^6  3^4  5^2  7^5  43 a^3  b^2  c  -2^8  5^4  7^7  a^5  e  \\
&  - 2^2  3^6  5^5  7^4  11  a^2  b  c  d +2^7  3^3  5^4  7^4  a^2  b^2  e+  2^3  3^3  5^4  7^5  491   a^4  c^2  -2^6  3^8  5^2  7^3  a  b^3  d  -2^5  3^9  5^3  7^2  b^3  f -2^3  3^5  5^6  7^3  a  b  d  e  \\
&   -2^2  3^6  5^7  7^2  b  e  f -2^{11}  3^3  7^7  a^5  b^2+ 2^2  3^2  5^7  7^2  601 a  c^2  e-2^{10}  3^2  5^2  7^8  a^6  c     +3^5  5^9  19  c^4+   2^2  3^9  5^7  7^2  17   c  d  f  -2\cdot 5^9  7^2  c  e^2  \\
\end{split}
\end{equation*}
\end{Small}
The invariant $I_{32}$ is an equation of degree 14, 8, 8, 4, 5, 4 in $a, b, c, d, e, f$ respectively.  We denote it as $I_{32}=\sum_{i=0}^{14} b_i J_2$ and display its coefficients as follows:
\begin{Small}
\[
\begin{split}
b_{14} = & -  2^{24} \cdot 3 \cdot 7^{18} \cdot c    \\
b_{13} = &   2^{21} \cdot 5 \cdot 7^{17} \cdot 41 \cdot e  \\
b_{12} = &   2^{20} \cdot 3^4 \cdot 5^2 \cdot 7^{16} \cdot c^2  \\
b_{11} = &  2^{16} \cdot 3 \cdot 7^{15} \cdot   c \left( 13824  b^2-140125 e \right)  \\
b_{10} = &   2^{14} \cdot 5 \cdot 7^{14} \cdot \left( -466560 bcd-874800 df-283392  b^2e-246375  c^3+401750  e^2 \right) \\
b_9 = &   2^{14} \cdot 3^2 \cdot 5^2 \cdot 7^{12} \left( 370440 dbe-90720 bcf-157248  b^2 c^2-510300 c d^2+1055875 e c^2  +595350  f^2  \right) \\
b_8 = &  2^{11} \cdot 3^3 \cdot 7^{10}  \left(  988722000  b^2ce+1190700000 b c^2d+3051168750 cdf+1190700000  d^2e -48771072  b^4c-1189015625  e^2c-100453125  c^4 \right)  \\
\end{split}
\]
\end{Small}

\noindent Next we describe $I_{40}$.  It is a degree 17 polynomial in $J_2$ and we denote it by $I_{40}=\sum_{i_0}^{17} A_i J_2^i$.  Then, we have

\begin{Small}
\begin{equation*}
\begin{split}
A_{17} & =   2^{33} \cdot 7^{22} \\
A_{16} & =  2^{28} \cdot 3^2  \cdot  5 \cdot  7^{21}  \cdot c^2  \\
A_{15} & =  2^{28} \cdot 3 \cdot 5^4 \cdot 7^{20} \cdot c e  \\
A_{14} & =  2^{22}\cdot 5 \cdot 7^{18}    \left(193536 b^2 e - 223425 c^3 - 266875 e^2  \right)   \\
A_{13} & =   2^{24}  \cdot 3^4 \cdot 5^2 \cdot 7^{17}  \left( - 168 b^2 c^2  - 1975 c^2 e - 1890  f^2  - 1680 b d e + 504 b {\it  c} f   \right)    \\
A_{12} & =     2^{18}  \cdot 3 \cdot 5^3 \cdot 7^{16}  \left( 4898880 c d  f    - 6531840 d^2 e + 2668750 c e^2  +  924075 c^4 + 326592 b c^2 d   -  1935360 b^2 c e  \right)    \\
A_{11} & = 2^{19}  \cdot   5  \cdot 7^{15} 
\left(  1020600000 b c d e + 96519600  b^2 c^3 - 172226250 c^2 d^2 - 41803776 b^4 e + 143184375 c^3 e  \right. \\
&   \left.      + 115290000  b^2 e^2 + 602791875 f^2 c - 71225000 e^3 - 306180000 b c^2 f + 1262992500 d e f \right)    \\
\end{split}
\end{equation*}
\end{Small}

\noindent The other coefficients are displayed in the Appendix~\ref{app}.

\begin{prob}
Express all invariants $I_8, I_{16}, I_{24}, I_{32}, I_{40}$ in terms of the transvectants of the binary octavics. 
\end{prob}

We summerize the above in the following theorem.

\begin{thm}   The  invariants  $J_2, \dots , J_8$ satisfy  the following equation
\begin{equation}\label{shaska}
 J_8^5 + \frac { I_8} {3^4 \cdot 5^3 } J_8^4 + 2 \cdot \frac { I_{16} } {3^8\cdot 5^6}   J_8^3 + \frac {I_{24}} {2 \cdot 3^{12} \cdot 5^6}  J_8^2  + \frac { I_{32}} {3^{16} \cdot 5^{10} }  J_8 + \frac { I_{40}} {2^2 \cdot 3^{20} \cdot 5^{12}}  =0, 
 \end{equation}
%
%where  $I_{32}$ and $I_{40}$ are displayed in the Appendix~\ref{app} and $I_8, I_{16}, I_{24}$ are as below. 
\end{thm}

\proof 
To prove that this relation holds we take a generic octavic \[f(x, z)= \sum_{i=1}^8  a_i x^i z^{8-i}.\] Compute the invariants $J_2, \dots , J_8$ and substitute them in the Eq.~ \eqref{shaska}.  We see that the equation is satisfied. This completes the proof. 

\qed

\begin{rem}
i)  In terms of the coefficients of the binary octavic the above equation is a degree 40 homogenous equation.

ii)  The equation    has   degrees in $J_i$, $i=2, \dots , 8$,  respectively    17, 10, 10, 6, 6, 5, 5.

iii) Similar relations as that in previous Theorem can be determined for $J_9$ and $J_{10}$ in terms of $J_2, \dots , J_7$.  However, such relations, as expected, are very large to display. 

iv)  In \cite{Shi1} it is commented that the field of fractions of $\cR$ is determined by a degree 5 equation 
\[ J_8^5 + a_1 J_8^4 + \dots \a_5 =0,\]
where $a_1, \dots , a_5$ are homogenous elements in $\Q [ J_2, \dots , J_7]$ but are not computed; see page 1043.  That equation is precisely Eq.~\eqref{shaska}. 

v) All coefficients of these equation can be expressed in terms of the transvectants of binary octavics. 

vi) The reader can check the correctness of the above equation in \cite{homepage} 
\end{rem}

% ************************
\begin{lem} The following hold true for any octavic. 

i)   If an octavic has a root of multiplicity exactly four then 
\begin{equation}\label{J_i}
\begin{split}
I_8 &= 2^{11}\cdot 3^6\cdot 5^4 \,  \cdot r^8,    
\quad       I_{16}=  2^{22}\cdot 3^{12} \cdot 5^7 \,   r^{16},      
\quad I_{24} =  2^{35}\cdot 3^{18} \cdot 5^7 \,    r^{24},   \\  
 I_{32} & = 2^{44}\cdot 3^{24} \cdot 5^{11} \,   r^{32}, 
\quad I_{40}=   2^{57}\cdot 3^{30} \cdot 5^{12} \,   r^{40}, 
\end{split}
\end{equation}
for some  $ r \neq 0 $.  Moreover, if the octavic has equation 
\[ f(x, y)=x^4 (a x^4+b x^3 y+c x^2 y^2+d x y^3+e y^4), \]
then  $r=e$. 

ii)  If an octavic has a root of multiplicity 5 then
\[ I_{8i}=0, \,  for \ \ i=1, \dots , 5.\]

\end{lem}

\proof i)   The proof follows Theorem~\ref{thm_5} part i) or by direct computation.  

ii)  Since all $I_{8i}$ for $i=1, \dots , 5$ are all homogenous polynomials in terms of $J_2, \dots , J_7$ then this is an immediate consequence of Theorem~\ref{thm_5}. 
 \qed

\begin{cor}
If $J_4=J_5=J_6=J_7=0$, then  $I_{24}=I_{32}=I_{40}=0$.  In this case, the Eq.~\eqref{shaska} becomes \[ J_8 \, \left( -10125 J_8+1075648 J_2^4-42336 J_3^2\, J_2  \right) =0\]
\end{cor}

The equation Eq.~\eqref{shaska} corrects the result of \cite[Thm.~5]{Shi1}.  To compute the other syzygies we follow a similar technique replacing $J_8$ by $J_9$ or $J_{10}$.  Indeed, both such cases are a bit easier from the computational point of view.  For our purposes of determining the field of invariants $\S_8$ the Eq.~\eqref{shaska} is enough.

\begin{figure}[htbp]
\begin{center}
\[
\xymatrix{ 
k[A_0, \dots A_8] \ar@{-}[dd]_{SL_2 (k)}   & \\
  & \S_8 = F(\cR_8) \ar@{-}[d]^5 \\
\cR_8 =P[J_8]  \ar@{-}[d]^5  \ar@{-}[ur] & F(P) \\
P=k[J_2, \dots , J_7] \ar@{-}[ur]&  \\
%k \ar@{-}[u]  &  
}
\]
\caption{The ring of invariants $\cR_8$ and its field of fractions}
\label{fig1}
\end{center}
\end{figure}
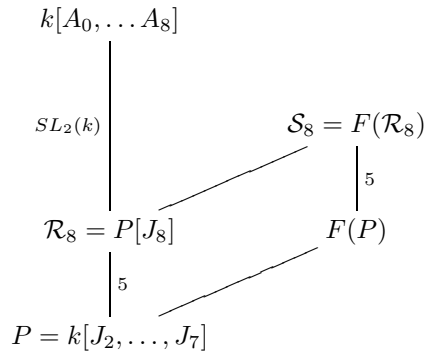

\begin{exa}
Let $C$ be the generic genus 3 hyperelliptic curve with automorphism group $\Aut(C) \iso \Z_2 \times D_8$.   Then $C$ has equation \[ C: \, \, y^2= x^8+\lambda \, x^4+1\]
Then,  invariants $I_{8i}$, $i=1, \dots 5$ of the corresponding binary form are given below.  
\[
\begin{split}
I_8 = & - 2^{11} \cdot 3^4 \cdot 5^4 \,     \left( 9\,{\lambda}^{2}-2450 \right)  \left( \lambda-14 \right) ^{3}
 \left( \lambda+14 \right) ^{3}
       \\
I_{16} = & 2^{22} \cdot 3^9 \cdot 5^7 \,   \left( 9\,{\lambda}^{2}+980 \right)  \left( 3\,{\lambda}^{2}-1960
 \right)  \left( \lambda-14 \right) ^{6} \left( \lambda+14 \right) ^{6
}
        \\
I_{24} = & - 2^{35} \cdot 3^{12} \cdot 5^7 \,      \left( 9\,{\lambda}^{2}-9310 \right)  \left( 9\,{\lambda}^{2}+980
 \right) ^{2} \left( \lambda-14 \right) ^{9} \left( \lambda+14
 \right) ^{9}
       \\
I_{32} = & 2^{44} \cdot 3^{16} \cdot 5^{11} \, \left( 9\,{\lambda}^{2}-12740 \right)  \left( 9\,{\lambda}^{2}+980
 \right) ^{3} \left( \lambda-14 \right) ^{12} \left( \lambda+14
 \right) ^{12}
          \\
I_{40} = & - 2^{57} \cdot 3^{21} \cdot 5^{12} \,   \left( 3\,{\lambda}^{2}-5390 \right)  \left( 9\,{\lambda}^{2}+980
 \right) ^{4} \left( \lambda-14 \right) ^{15} \left( \lambda+14
 \right) ^{15}
     \\
\end{split}
\]
If $\lambda^2=14$ then $I_{8i}=0$, for $i=1, \dots , 5$ dhe $J_i=0$ for $i=4, \dots , 8$.  This corresponds to the single curve with automorphism group $\Aut(C) \iso \Z_2\times S_4$.
\end{exa}

\begin{thm} Two genus 3 hyperelliptic curves $C$ and $C^\prime$ in Weierstrass form, given by equations  
\[  C:   Z^2=f(X, Y)   \textit{    and   } C^\prime:  z^2=g(X, Y) \]
are isomorphic over $k$ if and only if there exists some $\l \in k\setminus \{ 0\}$ such that 
\[ J_i (C) = \l^i J_i(C^\prime), \textit{   for   }  \,\,  i=2, \dots , 7,   \]  
and $J_2, \dots J_8$ satisfy the Eq.~\eqref{shaska}.
Moreover, the automorphism is given by
\[ 
\begin{aligned}
C & \to C^\prime \\
\begin{bmatrix} X \\ Y \end{bmatrix} & \to M \cdot \begin{bmatrix} X \\ Y \end{bmatrix}
\end{aligned}
\]
where $M  \in GL_2(k) $ and $\l = \left( \det M \right)^4$.
\end{thm}

\proof
The proof follows directly from the properties of invariants and Lemma~\ref{lem_3}. 

\qed

The above theorem gives a necessary and suffiecient condition for two hyperelliptic curves to be isomorphic.  However, $GL(2, k)$-invariants are prefered for identifying the isomorphism classes of curves.  In order to find such invariants we need to determine the field of fractions of $\cR_8=k[J_2, \dots , J_7, J_8]$.   

\subsection{On the invariant  field  of $GL_2(k)$}
%&&&&&&&&&&&&&&&&&&&&&&&&&&&&&&&&&&&&&&&&&&&&&&&&&&&&&&&&&

Let us assume that $J_2, J_3, J_4, J_5$ are all nonzero.  Define  the invariants 
\begin{equation}\label{abs_inv}
t_1:= \frac {J_3^2} {J_2^3}, \quad 
t_2:= \frac {J_4} {J_2^2}, \quad 
t_3:= \frac {J_5} {J_2\cdot J_3}, \quad 
t_4:= \frac {J_6} {J_2\cdot  J_4}, \quad 
t_5:= \frac {J_7} {J_2 \cdot J_5}, \quad 
t_6:= \frac {J_8} {J_2^4}. 
\end{equation}
Such invariants have the same degree in numerator and denumerator, therefore they are $GL(2, k)$-invariants.  Hence, $t_1, \dots , t_6 \in \S_8$.
For analogy with the genus 2 case, we call them \textbf{absolute invariants}.  For any two isomorphic genus 3 hyperelliptic curves $C$ and $C^\prime$ we have $t_j(C)=t_j(C^\prime)$, for $j=1, \dots , 6$. We would prefer an if and only if statement.

By substituting in Eq.\eqref{shaska}  get an affine equation of the hyperelliptic moduli of genus 3 as 
\begin{equation} \label{Sh3}
T(t_1, \dots , t_6 )=0
\end{equation}
This is an algebraic variety of dimension 5. It has degrees in $t_1, \dots , t_6$ respectively as 5, 10, 6, 6, 5, 5 and it has 25 464 terms. We denote this variety by $\T3$.  The equation of $\T3$ is explicitely computed and very useful in the arithmetic of genus 3 curves. The reader can check it at \cite{homepage}.

Then we have the following theorem.

\begin{thm}
The field of invariants of binary octavics is $\S_8= k(t_1, \dots , t_6),$
where $t_1, \dots , t_6$ satisfy the equation~\eqref{Sh3}. 
\end{thm}

\proof The proof of the theorem follows directly from Thm.~?? in Shioda. However, since that is based on Thm.~5 which contains syzygies which are incorrect, we provide a direct proof for this result.  

We denote the roots of the octavic $f(X,Y)$ by $(\alpha_i, \beta_i)$. Every $g\in GL_2(k)$ which fixes $f(X,Y)$ permutes these roots. Thus there is an $S_8$ action on $\{ \alpha_0, \dots , \alpha_7\}$. The  fixed
field is the invariant field of $GL_2(k)$ which we denote by $\S_8$. We can fix $\alpha_5=0$, $\alpha_6=1$, and $\alpha_7=\infty$. Let $s_1, \dots , s_5$ denote the  symmetric polynomials of $\a_0, \dots , \a_4$. Then
$[k(\a_0, \dots , \a_4):k(s_1, \dots , s_5)]=120$. Thus, $[k(s_1, \dots , s_5) : \S_8]=6\cdot 7\cdot 8= 336$.
\[
\xymatrix{ k(\alpha_0, \dots , \alpha_4)  \ar@{-}[d]^{\, \, \, 120}  \\
k(s_1, \dots , s_5) \ar@{-}[d]^{\, \, \, 336}\\
\S_8 \ar@{-}[d]   \\
 k(t_1, \dots , t_6) \ar@{-}[d]^{\, \, \, 5} \\ 
k(t_1, \dots , t_5)
}
\]
Our goal is to determine $\S_8$.

Since $t_1, \dots , t_6$ are $GL(2, k)$-invariants then $k(t_1, \dots , t_6) \subset \S_8$.   We know that $[k(t_1, \dots , t_6):k(t_1, \dots , t_5)=5$, since the degree in $t_6$ of the irreducible polynomial from Eq.~\eqref{Sh3} is 5.  If we  show that 
$ [k (s_1, \dots , s_5):k(t_1, \dots , t_5)] =5 \cdot 336,$ or equivalently $[k( s_1, \dots , s_5):k(t_1, \dots , t_6)]=336$ 
 then we are done. 
 
 The proof is computational.  Compute $t_1, \dots , t_5$ in terms of $s_1, \dots , s_5$.  This is computationally easy and we do not display these expressions here.  
 By Bezout's theorem we know that the degree $d=[k (s_1, \dots , s_5):k(t_1, \dots , t_5)]  $  is $d \leq 6 \cdot 4 \cdot 5 \cdot  6 \cdot 7$, becaus ethe degrees of $i_1, \dots , i_5$ are respectively 6, 4, 5, 6, 7.    There is at least one more solution at infinity.  Moreover, $d$ must be divisible by $5\cdot 336$.  Hence, $d= 5\cdot 336$ or $d= 2\cdot 5 \cdot 336$.  
 
From the system of equations we eliminate first  $s_5$. Continuing via the resultants we eliminate also $s_1$ and $s_4$.  We are left with two equations of degree 36 and 56.  From Bezout's theorem, the degree $d \leq 36\cdot 56$  and divisible by 1680.  Hence $d=1680$ and the proof is complete.  

\qed

\begin{cor}
Two hyperelliptic genus 3 curves with nonzero  invariants $J_2, J_3, J_4, J_5$ are isomorphic if and only if they correspond to the same point on the algebraic variety $\T3$. 
\end{cor}

\proof The proof is an imediate consequence of the previous Theorem.  \qed

%For more details of how explicitely this theorem can be used and how it is implemented in a computational package see next section.  

%We have a mapping 
%
%\[ 
%\begin{split}
%\H_3 \to \T3  \\
%[C] \to (t_1, \dots , t_6) \textit{ such that } T(t_1, \dots , t_6)=0 \\
%\end{split}
%\]

Since the moduli space of hyperelliptic curves is a rational variety then $\T3$ must have a birational parametrization.  Finding such parametrization via an equation of this size is very difficult.

%***************************************************************************
\subsection{Cases when $t_1, \dots , t_6$ are not defined}

To describe the moduli points in cases when absolute invariants are not defined is not difficult.  In this case, one has to treat each case separately when any of the invariants $J_2, \dots J_5$ are zero.  

Indeed, we can define invariants depending of which of the invariants is nonzero.  If  $J_2 \neq 0$, then  we  define
\[ i_1= \frac {J_3^2} {J_2^3}, \quad  i_2= \frac {J_4}  {J_2^2}, \quad  i_3= \frac {J_5^2}  {J_2^5}, \quad  i_4= \frac {J_6}  {J_2^3}, \quad  i_5= \frac {J_7^2}  {J_2^7}, \quad i_6=\frac {J_8} {J_2^4} \]

If $J_2 =0$ then we pick the smallest degree invariant among $J_3, \dots , J_7$ which is not zero.  This is possible because if $J_2=\dots = J_7=0$, then from Lemma~\ref{lem_4}    the binary octavic has a double root, hence we don't have a genus 3 curve.  For example, if $J_3 \neq 0$, then  we define
\[ i_1= \frac {J_2^3}  {J_3^2 }, \quad  i_2= \frac {J_4}  {J_3 \cdot J_4}, \quad  i_3= \frac { J_8}  {J_3\cdot J_5}, \quad  i_4= \frac {J_4^3}  {J_3^4}, \quad  i_5= \frac {J_5^3 }  {J_3^5}, \quad \]
Such invariants have high degree in some cases and therefore are not suitable for computations.  Hence, we prefer invariants $t_1, \dots , t_6$ defined in Eq.~\eqref{abs_inv}.

In the next example we see what happens in the case when all $J_2=J_3=J_4=J_5=0$.  We see that we get a genus 0 curve in the hyperelliptic moduli $\H_3$. 

\def\t{\tau}

\begin{exa}
Let us assume that $J_2=J_3=J_4=J_5=0$.  In this case $I_8$ and $I_8=I_{16}=0$ and 
\[  I_{24}-2679687500000u^4, \quad I_{32}-3204948120117187500\, u^3 \,v^2, \quad I_{40}=-306653442317962646484375 \,u^2 \,v^4 \]
where $u, v \in k[a,b,c,d,e]$.  Moreover, \[ ( 25 I_{24}\, I_{40} - 2 I_{32}^2 ) ( 25 I_{24}\, I_{40} + 2 I_{32}^2 ) =0\]
In this case, the Eq.~\eqref{shaska} becomes
\[ 2125764\, J_8^5-343000000\, J_6^4 \, J_8^2+16206750000 \, J_6^3\, J_7^2\, J_8-191442234375\, J_6^2\, J_7^4 =0\]
By defining 
\[ \t_1:= \frac {J_7^6} {J_6^7}, \quad \t_2 = \frac  {J_8^3} {J_6^4}, \]
the Eq.~\eqref{shaska} becomes 
\begin{Small}
\[ 
\begin{split}
-1064211156161261718750000000000\, \t_1\,\t_2+ 40353607000000000000000000\,{\t_2}^{2}+4649919888623184000000\,{\t_2}^{4}\\
-9606056659007943744\,{\t_2}^{5}-750282026508000000000000\,{\t_2}^{3}+7016382605513364494808197021484375\,{\t_1}^{2}\\
-19786546042268119734375000000\, \t_1\,{\t_2}^{2} & =0  \\
\end{split}
\]
\end{Small}
This is a genus 0 curve that can be parametrized as follows
\[ \t_1=  2^{46}\,  3^{6} \, 5^{60} \, 7^{36}   \cdot \frac { \left(t+ 2^{11}\,  3^{9} \, 5^{18} \, 7^{12} \right)^2 } {t^5},   \qquad
\t_2 =\frac {2^4\, 5^6 \, 7^3 } { 3^{12}} \cdot    \frac {\left(t+ 2^{11}\,  3^{9} \, 5^{18} \, 7^{12} \right)^2 } {t^2}   \]
It is possible in this case to express the equation of the curve in terms of the parameter $t$.  
\end{exa}

%****************************************************************************
\subsection{A computational package for genus 3 hyperelliptic curves}   

All the computational results described in this paper are implemented in a Maple package which will be made freely available. This package among other things computes the following:  \\

i) Invariants $J_i$ for $i=1, \dots , 10$. Their formulas are given in terms of the coefficients of a generic octavic \[f(X, Y) = \sum_{i=1}^8 a_i \, X^i Y^{8-i}\] and can be evaluated on any given octavic. 

ii) Invariants $I_{5i}$, for $i=1, \dots , 5$. Their formulas are given in terms of $J_2, \dots , J_7$ and can be evaluated on any octavic.  

iii)   The equation~\eqref{shaska} in terms of the invariants $J_2, \dots , J_8$.  

iv) The equation~\eqref{Sh3} in terms of invariants $i_1, \dots , i_6$.   \\

Some  problems which we are further studying are:  finding a minimal model of a genus 3 curve over its minimal field of definition, determining an algorithm which determines when the field of moduli is a field of definition, and describing the loci of curves with fixed automorphism group in terms of invariants $t_1, \dots , t_6$. For these problems and other computational aspects of genus 3 hyperelliptic curves see \cite{sh_thompson}.

%*********************************************************************************************************
\bibliographystyle{amsplain}

%****************************************************************
\appendix 

\section*{The algebraic relations among invariants}\label{app}

\begin{tiny}
\[
\begin{split}
I_{24} & =  -96060295941120  a^2b^4c^2+ 4881286149304320  a^5b^2c^2-4065020080128 ab^6c  +309844863885312  a^4b^4c-7872354689826816  a^7b^2c            \\
& -34767360928125 ab^2c^4-5975295019312500  a^3c^2d^2+16642058118840000  a^5  cd^2    -19599203957760 ab^4d^2+995929919631360  a^4b^2d^2                       \\
& -83701957500000 bc^4d+16332669964800 b^5cd-62983068750000 ac^4e - 7507549546875 b^2c^3e+1950385122000000  a^3c^3e           \\
& +40351024369152000  a^7ce-105389409484800  a^3b^4e+2677671663206400  a^6b^2e -2000752070688000  a^3bd^3+632931468750000 ac^3d^2       \\
& -65637337500000 ad^2e^2-495116212500000  a^2c^2e^2+3237512241600000  a^4ce^2 +350688139200000  a^3b^2e^2+ 182284263000000 bd^3e            \\
& +5907570414480000  a^4d^2 e+265581311901696  a^4b^3f-3373866295640064  a^7bf +76081687500000 ace^3+18074283726643200  a^6df                      \\
& -5226454388736 ab^5f-190547816256000  a^2b^2f^2+164093343750000 de^2f  + 1435488571125000 ad^3f-452093906250000 c^3df                                     \\
& +174583596796875 ac^2f^2-116965735425000 b^2cf^2+1813953459150000  a^3c f^2   +136256080078125 cef^2-551353635000000  a^2ef^2                       \\
& -183784545000000  a^2bcef-1727902909687500 acdef+2507801378040000  a^3bcde  -623601590062500 abc^2de- 7082804531250 bc^2ef                    \\
& -74939085900000 b^2def-43758225000000 abe^2f +66672040953348096  a^{10}c  + 382796952300000 abdf^2+90329478750000 bcde^2                           \\
& -63515938752000    a^2 b^3 c f- 1506298130820000  a^3 b c^2 f-98703552780000 b^3 c d e  -232513704360000 a b^2 d^2 e-199865692687500  a^2 cd^2 e           \\
& +1613775332736000  a^5bcf  +3447798064200000  a^3def  +1327029434640000 ab^2cdf  + 51690608865000 b^3c^2f +150697968750000 c^2d^2e   \\
& +10543332173875200  a^6bcd-829941599692800  a^3b^3cd- 7284631253256000  a^4bc^2d +286712891892000 ab^3c^2d   +5769351562500 c^6     \\
& +478496190375000 abcd^3-3176319702528000   a^4b^2ce+62507749248000 ab^4ce +618356229120000  a^2b^2c^2e+363935594531250 abc^3f               \\
& -30474418113720000  a^4cdf+5709463802437500  a^2c^2df-535915733220000  a^2bd^2f -444317891062500 bcd^2f+768464444160000  a^4bef          \\
& -127423951200000 a b^2 c e^2  -22677564950118400  a^9e  -30245685120000 ab^3ef  - 1956290913561600  a^3b^2df   -655008118380000  a^2b^2cd^2  \\
& -2679687500000 e^4-4518527895000 b^4c^3+15303471326640000  a^6c^3-62010412933754880 {J2}^8c^2 -4455038212800000  a^6e^2   \\
& -12651998608650240  a^7d^2-72461536875000  a^2c^5-839418892200000  a^4c^4-247062900000000  a^3e^3 -2956078125000 c^3e^2      \\
& +78746801616000 b^3d^3  +4841325998208000  a^5f^2  -256337244843750 d^2f^2 +113927664375000 bf^3 -2009683999575000  a^2d^4  \\
& -487519136370000  a^3b^2c^3     +9724050000000 b^2e^3  +18391180106250 b^2c^2d^2   -15710828636160000  a^5c^2e +48998009894400 b^4df        \\
&  -6901297200000 b^4e^2 +1382659891200 b^6e   +1139135992312500  a^2bc^3d  \\ \\
b_7 = & 2^8\cdot 3\cdot 5\cdot 7^{10}  \left( 2821754880  b^3cd-23346225000 def+7031205000 b c^2f+5290790400 fd b^2   +469476000  b^2 c^3+19190925000  c^2 d^2 \right) \\
& \left. -7589671875  c^3e+856977408  b^4e-2429784000  b^2 e^2  -27938925000  f^2c+1599500000  e^3-20182365000 bcde  \right) \\
b_6 = & 2^6\cdot 3^5\cdot 5^2\cdot 7^8 \left( + 1590540000 db e^2 -1793022000  b^2 c^2e-995742720 d b^3e-8574693750 c d^2e-1574370000 b c^3d  -2601112500  c^2df \right. \\
&  -3429216000 b d^2f  +243855360  b^3cf-1371686400  b^2 f^2  +2976750000 e f^2+130056192  b^4 c^2+2192312500  e^2 c^2+276921875  c^5  \\
& \left. -1199520000 bcef  \right) \\
b_5 = & 2^5\cdot 3 \cdot 7^7  \left(91423434375000 b c^2de  -7688302272000  b^4ce-18517766400000  d^2 b^2e -18517766400000  b^3 c^2d -64531023046875  c^3 d^2 \right. \\
& -52081218000000 bc d^3+18491571000000  b^2c e^2-36095751562500 b c^3f +26040609000000  f^2db-97652283750000 f d^3\\
&  +2775093750000 e c^4+46632649218750  d^2 e^2-10697312500000  e^3c    +81089859375000  c^2 f^2+143497954687500 cdef   \\
& \left.  +252829237248  b^6c +3794809500000  b^2 c^4 -47451776400000  b^2cdf \right)\\
b_4 = & 2^3\cdot 5 \cdot 7^5 \left(2563321813920000  b^3cde -2551979682000000 b f^3+11278838773125000  d^2 f^2 +93467500000000  e^4 -335986353561600  b^4df \right. \\
& -36280995545088  b^6e+154301003136000  b^4 e^2-203149296000000  b^2 e^3-1549416235500000  d^3be -5506503778125000 d e^2f\\
& -63420152343750  c^6-673882453125000  c^3 e^2-1958571267187500 ce f^2+290775637687500  b^2 c^3e -179192721899520  b^5cd \\
& +339070429687500  c^3df-893019284640000  b^3 c^2f+3183898460400000  b^2c f^2+34998578496000  b^4 c^3 +2965157344800000  b^2def \\
&    -3076046938125000 c d^4+2049090513750000 b c^2ef+11757334963500000 bc d^2f -249989846400000  b^2 c^2 d^2\\
& \left. +6684585148828125  c^2 d^2e +253301757187500 b c^4d -3137933570625000 bcd e^2  \right)    \\
b_3 = & 2^2\cdot 3^3\cdot 5^2\cdot 7^5   
\left(-119352791250000  f^3d+923835937500  c^5e+2036562500000  c^2 e^3 -14145516000000  b^2e f^2+1453824288000  b^4 c^2e\right. \\
& +1935723847680 d b^5e+2263282560000  b^3 c^3d  +2666558361600  b^4c d^2+42160986000000 b c^2 d^3+215648793281250 c d^3f\\
& -6184019520000 d b^3 e^2-11786069531250 c d^2 e^2+3942540000000 db e^3-474054819840  b^5cf  +7343991562500 b c^4f-77308057968750  d^4e \\
& +2498818359375  c^4 d^2-13526936718750  c^3 f^2+11627929687500  e^2 f^2-19656915468750 b c^3de  +11107444950000  b^2c d^2e\\
& +12296954250000 bcd f^2 -32049232031250  c^2def-55456852500000 b d^2ef+12579405300000  b^2 c^2df+2222131968000  f^2 b^4\\
& \left. -4140860062500  b^2 c^2 e^2 -63207309312  b^6 c^2-807132093750  b^2 c^5-955040625000 bc e^2f+4663733760000  b^3cef+13332791808000 f d^2 b^3    \right) \\
b_2 = &  2\cdot 3 \cdot 7^3  
\left( 13999431398400000  b^4 d^2e-501993771152343750  d^3ef+147650253030000000  b^2 d^3f +13999431398400000  b^5 c^2d\right. \\
& -70508565618750000  b^2 d^2 e^2-39373400808000000  b^3d f^2+640843112109375000  f^4 -113927664375000000 bde f^2 +199373412656250000 bc f^3  \\
& -219625488281250  c^7-1121475446191406250 c d^2 f^2-113792664111328125 e c^2 f^2 +3874904345088000  b^6ce+11248855857421875  b^2 c^4e\\
& -4130738964843750 b c^5d+12125158565625000  b^2 c^3 d^2+78746801616000000  b^3c d^3 -13979627676000000  b^4c e^2+16174336500000000  b^2c e^3\\
& -95569451679744  b^8c-4556693232000000  b^4 c^4+357041162460937500  c^2 d^4+5361029296875000  c^4 e^2 -204478306347656250  c^3 d^2e \\
& -7851484375000000 c e^4-216968907487500000  b^2cdef-138232232775000000  b^3 c^2de -85364120214843750 b c^3ef+35873542958400000  b^4cdf \\
& +54576776362500000  b^3 c^3f +81255806542968750  c^4df -108366909328125000  b^2 c^2 f^2  -93278275207031250 bc d^3e  \\
& \left.  -372604495183593750 b c^2 d^2f +146072378320312500 b c^2d e^2+40730312109375000  d^2 e^3 +377268178710937500 cd e^2f \right) \\
\end{split}
\]

\[
\begin{split}
b_1 = &  5\cdot 7^2    \left(  5359375000000000  e^5+3951865858007812500 bc d^2ef+979586879717376  b^8e  +275613805656000000  b^3 f^3-9768884765625000  c^6e\right.   \\
&    +77017426171875000  c^5 d^2    -5554836112896000  b^6 e^2+2321709802734375  b^2 c^6+10970061984000000  b^4 e^3   +4103367187500000  c^3 e^3\\
& +55159493115234375  c^4 f^2+167336953434000000  b^3 d^3e    +29607242411250000  b^3 c^4d-104620750718400000  b^4 c^2 d^2-1021687018734375000 b c^3 d^3  \\
& +231488824218750  b^2 c^3 e^2-229060912500000000  c^2 d^2 e^2   +48223041370560000  b^5  c^2f+11649491191406250 b c^5f+12095508728217600  b^6df    \\
& +6450937988382720  b^7cd -6756317381953125000  c^2 d^3f +22785532875000000 b d^3 e^2   -152243816457600000  b^4   c f^2   -3593187392256000  b^6 c^3  \\
& -885901518180000000  b^2 d^2 f^2+5415124297324218750  d^2e f^2-157134027832031250 c e^2 f^2   +797493650625000000 be f^3+20649479586750000  b^4 c^3e \\
& +338896825627500000  b^3cd e^2-179760351093750000 bcd e^3-138419377951680000  b^5cde  -383861214843750000 d e^3 f  -865138201347656250 cd f^3  \\
& +91935556305468750  b^2 c^2 d^2e-3249990068554687500 b c^2d f^2+325453361231250000  b^2ce f^2     +155009007070312500 b c^4de +594702408037500000  b^2d e^2f  \\
& -10094490000000000  b^2 e^4  -476970373441406250  b^2 c^3df -1269792176058000000  b^3cd^2f  +632366351367187500  c^3def \\
& \left.  +3076046938125000000 c d^4e+12118628320312500 b c^2 e^2f-221301775485000000  b^3 c^2ef-160118496619200000  b^4def  \right)   \\
b_0 = &  3^2\cdot 5^2    \left( 27459204328828125  b^3 c^3de+44149248498600000  b^4c d^2e+106744394531250 b c^4ef +7188851250000000 b c^3d e^2\right. \\
&  +440216495145000000  b^3 d^2ef-1196240475937500000 c d^3ef-59716324558800000  b^4 c^2df -516775885605000000  b^2c d^3f-375961292437500000 b c^2 d^3e\\
& -256970176312500000 b d^2 e^2f-379845703125000 bc e^3f-183742537104000000  b^3 cd f^2 +961264668164062500 b c^3 d^2f -131734136362500000  b^2c d^2 e^2  \\
& +7581112481250000  b^3c e^2f-18510359293440000  b^5cef+344517257070000000  b^2d f^3 -52917850685952000  b^5d^2f-5121925300961280  b^7de\\
& +27800635693359375  c^3e f^2+2210908736777343750  c^2 d^2 f^2-62765703984375000  b^2 e^2 f^2 +19967480859375000  c^4 d^2e +344517257070000000  b^2 d^4e  \\
& +66351471732000000  b^4e f^2-59902442578125000  c^5df-11549922890625000 b c^6d  -14111426849587200  b^6c d^2 -162415278333000000  b^3 c^2 d^3  \\
& -1897468822265625  b^2 c^5e+24544373474880000  b^5d e^2-960461782875000  b^4 c^2 e^2 +1529351464032000  b^6 c^2e +21879112500000000 c d^2 e^3 \\
& -31295882520000000  b^3d e^3+1856179335937500  b^2 c^2 e^3+4456856250000000 bd e^4 +1254349053296640  b^7cf+36053358788671875  b^2 c^3 f^2\\
& -40051757812500  c^5 e^2-334960937500000  c^2 e^4-186912574365234375 c f^4 +6328909611360000  b^5 c^3d +287652282750000  b^4 c^5 +14244213867187500  e^3 f^2  \\
& -809954488916015625 de f^3 -24614001562500000  c^2d e^2f -4703808949862400  b^6 f^2+379628173828125  c^8-284937835587890625 bcde f^2 \\
& \left.  -4141406228203125  b^3 c^4f +92505194627343750  b^2 c^4 d^2+65342262505078125  b^2 c^2def  -33110227458984375 b c^2 f^3 -83623270219776  b^8 c^2\right)  \\
A_{10 } & =    2^{14}  \cdot 3^2 \cdot 5^2 \cdot 7^{14} \ 
\left( 682668000 b c e f + 41803776 b^4c^2 - 1131034375 c^2e^2 + 940584960 f^2b^2 - 2483460000 f^2e - 200559375 c^5 - 352836000 b c^3d   \right. \\
& \left. + 680486400 b^2c^2e + 836075520 b^3d e + 7144200000 c d^2e - 1152900000 b d e^2 - 250822656 b^3c f - 5292540000 c^2d f   \right)     \\
A_{9 } & =     2^{14}  \cdot 3^2 \cdot 5^3 \cdot 7^{12} \cdot  \left (- 38623331250 c d e f - 4938071040  b^2c d f - 14979006000 b c^2d e - 1332223200 b^2c^4 + 6573301875 c^3d2 + 851175000 c^4 e  \right. \\
&  - 11348859375 d^2e^2 + 2492875000 c e^3 + 6697687500 f^2c^2 - 658409472  b^3c^2d + 1463132160 b^4c e + 4938071040 b^2d^2e - 4035150000 b^2c e^2 \\
& \left. + 3597615000 b c^3f  - 9258883200 f^2d b    \right)    \\
A_{8 } & =     2^{10}   \cdot 5 \cdot 7^{10}  \left( 2624893387200000 b c d^2f   - 846721653750000 b c^2e f - 1283281211520000 b^2d e f  + 1429956281250000 b c d e^2  + 311098475520000  b^3c^2f \right. \\
& - 1036994918400000 b^3c d e - 49035045427200 b^4c^3 + 14158437285888  b^6e - 58571009280000 b^4e^2 + 184454313750000 b c^4d + 174992892480000 b^2c^2d^2 \\
& + 72369158400000  b^2e^3 - 1640558367000000 d^2f^2 - 145827410400000 b f^3 + 19812568359375 c^6    - 20460988800000 b^2c^3e - 3106789323750000 c^2d^2e  \\
&  + 229218018750000 c^3e^2 - 29292200000000 e^4 - 2916548208000000 b d^3e   + 1515418773750000 c^3d f + 1476840093750000 d e^2f - 612475123680000 b^2c f^2\\
&  \left.   + 2253236028750000 c e f^2   \right)    \\
A_7 & =  2^9  \cdot 3^2 \cdot 5^2 \cdot 7^{10} 
\left(  78121827000000 f^3 d - 4266493378560 f^2  b^4 - 4431306250000 c^2 e^3 - 26040609000000 d^4e + 1137731567616 b^5 c f\right. \\
& - 1798832812500 c^5 e - 17034451875000 f^2e^2 - 29206701562500 f^2c^3 + 34342392656250 c d^2e^2   - 6461532000000 b d e^3 + 3200928192000 b^3c^3d  \\
& - 3792438558720 d  b^5 e  + 10459108800000 d b^3e^2 + 6477790725000 b^2c^2e^2 + 10725230250000 b c^3d e  - 18517766400000 b^2c d^2e - 126414618624 b^6c^2  \\
&  - 1714999910400 b^4 c^2 e+ 793334250000 b c^4 f + 22529949120000 b^2 e f^2 - 6515380546875 c^4d^2 - 20832487200000 b c^2d^3 + 884987775000 b^2c^5 \\
& \left.  - 27776649600000 b c d f^2 + 3819170250000 b c e^2 f - 6193164096000 b^c e f + 61900984687500 c^2 d  e  f   + 114578679600000 b d^2 e f + 20237273280000 b^2c^2d f    \right)    \\
A_6 & =  2^9 \cdot 3 \cdot 5^3 \cdot 7^8 \cdot 
\left( 853442578125 c^7 - 546852789000000 f^4 + 251989765171200 f^2 d b^3 - 1806131250000 c^4 e^2 + 18307625000000 c e^4  \right. \\
& - 476054883281250 f^2 e c^2+ 259612414031250 b c^3j 6 f - 1781177655600000 b c^2 d^2 f + 2187411156000000 b c d^3 e - 433879545093750 b c^2 d e^2 \\
&  + 1015930959120000 b^2 c d e f+ 245638171296000 b^3 c^2 d e + 10906964409600 b^4 c^4 - 461407040718750 c^2 d^4 - 147215171250000 d^2 e^3 \\
& + 298984770000000 c^3 d^2 e- 285925886820000 b^2 c^2 f^2 - 1589290918031250 c d^2 f^2 + 856736036100000 b c f^3 + 8959636094976  b^5 c^2 d  - 45230724000000 b^2 c e^3\\
& + 36606880800000 b^4 c e^2- 557331321093750 c d e^2 f - 37629282796875 b c^5 d - 109370557800000 b^2 c^3 d^2 - 8849023303680 b^6 c e - 34426624050000  b^2 c^4 e \\
& \left. + 44798180474880 b^4  d^2 e- 49303553040000 b^3 c^3 f + 62292081796875 c^4 d f + 3349473332625000 d^3 e f - 1130162430600000 b d e f^2  \right)    \\
A_{5 } & =     2^5 \cdot 5 \cdot 7^7  
\left(1417442429088000000 f^2d^2b^2  - 12485200000000000 e^5  + 22840481250000000 c^6e + 83710757812500000 c^3e^3 - 31263476428800000 b^4e^3\right. \\
&  + 125994882585600000  f^3b^3 + 843381182109375000 f^2e^2c + 470581317773437500 f^2c^4  - 7362005671912500000 b c d^2e  f + 3166863562012500000 b c^2d f^2  \\
& - 1946795928840000000 b^2c e f^2 - 451139485781250000 b c^2e^2f - 1275989841000000000 b^2d e^2f + 731567508840000000 b^3c^2e f + 554377483376640000 b^4d e f  \\
& - 865759830468750000 c^3d e f - 2267907886540800000 b^3c d^2f + 331236546480000000 b^2c^3d f - 1235482227000000000 b^3c d e^2 + 763268467500000000 b c d e^3 \\
& + 16868450672640000 b^6e^2 + 83335976718750000 b c^4d e - 49997969280000000 b^2c^2d^2e + 447981804748800000 b^5c d e - 159368527080000000 b^3c^4d \\
& - 75596929551360000 b^4c^2d^2 + 2627822955712500000 b c^3d^3 - 45171004838400000 b^4c^3e + 2519897651712000000 b^3d^3e + 1538023469062500000 c d^4e\\
& - 3058222453751808 b^8e- 134394541424640000 b^5c^2f - 201084911015625000 b c^5f + 6196037403937500000 c^2d^3f + 264589253429760000 b^4c f^2 \\
& + 8688968085937500  b^2c^6+ 894543142500000000 d e^3f - 22660212444187500000 d^2e f^2 + 8920536120562500000 c d f^3 + 364568526000000000 b e f^3 \\
&  - 49090588903125000 b^2c^3e^2 - 932669728593750000 c^2d^2e^2 - 3474793763437500000 b d^3e^2+ 14122093083033600 b^6c^3 - 262331525654296875 c^5d^2\\
& \left.  + 25308460800000000 b^2e^4    \right)    \\
\end{split}
\]

\[
\begin{split}
A_4 & =      2^4  \cdot 3^4 \cdot 5^2 \cdot 7^5 \cdot 
\left( 6689861617582080 f^2b^6  - 231356250000000 c^8 - 428732586576000000 f^3d b^2  - 3309682031250000 c^5e^2    \right. \\
& - 3433215625000000 c^2e^4 - 358872142781250000 f^4c + 303638323275000000 b c d e f^2 - 451862434170000000  b^2c^2d e f + 148663591501824 b^8c^2 \\
& + 71562930468750000 c d^2e^3+ 116421883200000 b^4c^5 + 244588454455078125 c^3d^4 - 422283964306640625 d^4e^2 + 130661359718400000  b^3c d f^2 \\
& + 81695043281250000 f^2e c^3- 17965376856000000  b^3c e^2f - 4063263750000000 c^2d e^2f + 620272839375000000 b d^2e^2f + 1415389500000000 b  c e^3f \\
& - 12302724000000000 b d e^4 - 43197572742187500 c  d^3e f + 367485074208000000 b^2c  d^3f + 14566321953792000 b^5c  e f - 38516912437500000 b c^4e  f \\
& - 538978108838400000 b^3d^2e  f + 402924173006250000 b c^3d^2 f + 17732613104640000 b^4c^2d  f - 818506975387500000 b c^2d^3e \\
&  - 1783963098021888 b^7c f+ 73723926937500000 b c^3d e^2 + 138744204007500000 b^2c d^2e^2 + 2058386904000000 b^3c^3d e - 65330679859200000 b^4c d^2e \\
&  - 23904956250000000 f^2e^3 + 538315584307200 b^6c^2e + 5861394393750000 b^2c^5e + 5946543660072960 b^7d e + 10913044570312500 c^4d^2 e \\
& + 30395046528000000 b^3d e^3 - 183742537104000000 b^2d^4e - 6338258035200000 b^4c^2e^2 - 24599823897600000 b^5d e^2 + 9851285850000000 b^2c^2e^3 \\
& - 19484805312000000 b^3c^4f + 6963501972656250 c^5d   f + 204183691818750000 b^2c^3f^2 - 1342922343820312500 c^2d^2f^2 \\
& - 52990440330240000 b^4e f^2 + 80130061620000000 b^2e^2f^2 - 444824236237500000 b c^2f^3 + 1913984761500000000 d e f^3 \\
& \left. - 7528583107584000 b^5c^3d + 6870810761718750 b c^6d - 10365609082500000 b^2c^4d^2 + 97996019788800000  b^3c^2d^3  \right)   \\
A_3 & =     2^3  \cdot 3 \cdot 5^3 \cdot 7^5 \cdot 
\left (  - 231356250000000 c^8 - 428732586576000000 f^3d b^2 + 6689861617582080 f^2b^6  - 3309682031250000 c^5e^2 + 97996019788800000  b^3c^2d^3  \right.  \\
& - 3433215625000000 c^2e^4- 358872142781250000 f^4c + 303638323275000000 b c d e f^2 - 451862434170000000  b^2c^2d e f + 148663591501824 b^8c^2 \\
&  + 81695043281250000 f^2e c^3 + 116421883200000 b^4c^5 + 244588454455078125 c^3d^4 - 422283964306640625 d^4e^2 + 130661359718400000  b^3c d f^2 \\
& + 30395046528000000 b^3d e^3 - 17965376856000000  b^3c e^2f - 4063263750000000 c^2d e^2f + 620272839375000000 b d^2e^2f + 1415389500000000 b c e^3f \\
& + 71562930468750000 c d^2e^3 - 43197572742187500 c  d^3e f + 367485074208000000 b^2c  d^3f + 14566321953792000 b^5c  e f - 38516912437500000 b c^4e  f \\
& - 538978108838400000 b^3d^2e  f + 402924173006250000 b c^3d^2 f + 17732613104640000 b^4c^2d  f - 818506975387500000 b c^2d^3e  - 7528583107584000 b^5c^3d \\
&  - 12302724000000000 b d e^4 + 73723926937500000 b c^3d e^2 + 138744204007500000 b^2c d^2e^2 + 2058386904000000 b^3c^3d e - 65330679859200000 b^4c d^2e \\
&  - 23904956250000000 f^2e^3 + 538315584307200 b^6c^2e + 5861394393750000 b^2c^5e + 5946543660072960 b^7d e + 10913044570312500 c^4d^2 e \\
& - 183742537104000000 b^2d^4e - 6338258035200000 b^4c^2e^2 - 24599823897600000 b^5d e^2 + 9851285850000000 b^2c^2e^3 + 6870810761718750 b c^6d\\
& - 1783963098021888 b^7c f- 19484805312000000 b^3c^4f + 6963501972656250 c^5d  f + 204183691818750000 b^2c^3f^2 - 1342922343820312500 c^2d^2f^2 \\
& \left. - 52990440330240000 b^4e f^2 + 80130061620000000 b^2e^2f^2 - 444824236237500000 b c^2f^3 + 1913984761500000000 d e f^3 - 10365609082500000 b^2c^4d^2 \right) \\
%&        \\
%
A_2 & =  2^3 \cdot 7^3 \left( 64873168945312500 c^9 - 18757812500000000000 e^6 - 1339483095488160000000 f^2d^2b^4 + 26161779208753125000000 f^4d b \right.   \\
&  - 119065164043392000000 f^3b^5 + 49053336016412109375000 f^3d^3  + 11537813232421875000 f^2e^3c  - 560843682604980468750 f^2e c^4 \\
&  + 29542831787109375000 e^2c^6 - 10292009765625000000 e^4c^3  + 1910513527976074218750 f^4c^2  + 13914190719914625000000 b^3c d^2 e f \\
& + 6991331335149902343750 b^2c^3d e f + 10034001027852539062500 b c^2d e f^2 - 22244465475206542968750 b c d^2e^2f - 6672688981733376000 b^8c^3 \\
& - 1593990434186718750000 b^2c e^2f^2 - 32598715496484375000 b^4c^6  + 1251360612663574218750 c^4d^4 + 1156008087518183424 b^{10}e + 98106672032824218750000 d^6e \\
& + 42827801519514375000000 b^2d^2e  f^2- 7970342942822400000 b^8e^2 + 19695990150144000000 b^6e^3 - 23916495456000000000 b^4e^4 + 23597028000000000000 b^2e^5  \\
& + 16569799717478027343750 c d  e f^3 - 67175712597656250000 b c^2e^3f - 1690686539325000000000 b^2d e^3f - 5985372132203625000000 b^3c^2d f^2 \\
& - 88296004829541796875000 b c d^3f^2 + 1839722152753800000000 b^4c e f^2 + 71944892824071093750000 b d^4e f + 852653628126562500000 b^3c^2e^2f  \\
& + 1205810399745000000000  b^4d e^2f - 1042359012597656250000 c^3d e^2f + 46248846347043457031250 c^2d^3e f - 29151696832610625000000 b^2c^2d^3f \\
& + 62484128808192000000 b^6c^3e- 691331295853800000000 b^5c^2e f - 91760935034179687500 b c^5e f - 349257814527283200000 b^6d e f - 1863346193238600000000 b^4c^3d f \\
& + 2143172952781056000000 b^5c d^2 f - 1497513029475585937500 b c^4d^2 f + 12301784590829589843750 b c^3d^3 e + 11627457426112500000000 b^2c d^4e \\
& + 1167530704515000000000 b^5c d e^2  - 361957683691406250000 b c^4d e^2 - 5363296110632373046875 b^2c^2 d^2e^2 - 1442577403575000000000 b^3c d e^3 \\
& + 583898814843750000000 b c d e^4 - 780546910549218750000 b^3c^4d e + 2631127508994600000000 b^4c^2d^2e - 282228536991744000000 b^7c d e \\
& + 443216854028320312500 c^2 d^2e^3 - 4057248947554687500000 b  d^3e^3 + 84668561097523200000 b^7c^2f + 566963056184765625000 b^3c^5f \\
& + 143385921661376953125 c^6d f - 147160008049236328125000 c d^5f - 166691229660748800000 b^6c f^2  + 1784515113281250000000 d e^4f \\
& - 2439437968006347656250 b^2c^4f^2 - 5677469446343994140625 c^3d^2f^2 - 49947712684749755859375 d^2e^2f^2 + 3532647655502929687500 b c^3f^3 \\
& - 996867063281250000000 b e^2f^3 - 689034514140000000000 b^3e f^3  + 150603258090600000000 b^5c^4d - 40302456207275390625 b c^7d \\
&  + 1258810413841845703125 b^2c^5d^2 - 4966585386296625000000 b^3c^3d^3 - 9810667203282421875000 b c^2d^5 + 47626065617356800000 b^6c^2d^2\\
& - 13952948911335000000000 b^2c d f^3  + 1746549092285156250 b^2c^6e - 685505045269775390625 c^5d^2e - 2381303280867840000000 b^5d^3e \\
& \left. - 94370714922093750000 b^4c^3e^2 + 6567360212896875000000 b^3d^3e^2 - 15476361157441406250000 c d^4e^2  + 15708519052734375000 b^2c^3e^3  \right)  \\
A_1 & =     2  \cdot 3^2 \cdot 5 \cdot 7^2 
\left (44253777130664062500\,b^2 c^2d\,e^2f  - 1475682251116500000000\,b^3c\,d\,e\, f^2 - 478911551651367187500\,b\,c\,d\, e^2 f^2    \right.\\
&  - 10916890583405625000000\,b^2c\,d^3 e\,f + 1488478606379100000000\,b^4c^2 d\,e\,f + 1190651640433920000000\,f^3d\,b^4 - 552727184194335937500\,c^3d\,f^3  \\
& - 21801482673960937500000\,f^4d^2 - 8720593069584375000000\,f^5b - 43602965347921875000000\,d^5e\,f  + 2054184145279875000000\,b^3c^2 f^3 \\
& + 65404448021882812500000\,c\, d^4 f^2 + 85844513334988800000\,b^6 e\,f^2 - 194716049736600000000\,b^4e^2 f^2  + 3457882625756835937500\,d\,e^2f^3 \\
&  - 2884862076345703125000\, b\,c^3d^2e\,f - 8128181865362227200\,f^2b^8 - 658479114554625000000\,b^4c^3f^2 - 193790957101875000000\,b^2c\,f^4 \\
&     - 4940454093750000000\,b^2c^2e^4   - 56323365820312500000\,c^4d^2e^2 + 2775771747000000000\,b^4c^2e^3 - 73859963063040000000\,b^5d\,e^3 \\
& + 116178087375000000000\,b^2e^3 f^2- 29496285000000000000\,b\,d\,e^5 - 62902448437500000000\,c\,d^2e^4 + 59791238640000000000\,b^3d\,e^4 \\
\end{split}
\end{equation*}

\begin{equation*}
\begin{split}
& + 2167515164096593920\,b^9c\,f + 85627772138160000000\,b^5c^4f - 2370465593261718750\,b\,c^7f - 1709677874091796875000\,c^4d^3f  \\
& - 262287930884765625000\,b\,c^5d^3  + 1959173759434752000\,b^8c^2e + 1305828038812500000\,b^4c^5e - 7225050546988646400\,b^9d\,e  \\
& - 65891601562500000\,e^3c^5 - 4145981297976000000\,b^6c^2e^2 + 769011734531250000000\,c^3d^4e + 3973403109375000000\,b^2c^5e^2  \\
& + 1190651640433920000000\,b^4d^4e+ 12196304634286080000\,b^7c^3d - 19151182255078125000\,b^3c^6d + 124852350865725000000\,b^4c^4d^2  \\
& - 6878792970000000000\, b^3c\,e^3f+ 1453432178264062500000\,b^2c^3d^4 - 6394981327382812500\,b^3c^3d\,e^2 - 1066855106060100000000\,b^4c\,d^2e^2  \\
& + 797493650625000000000\, d^4e^3 + 39851714714112000000\,b^7d\,e^2 - 238130328086784000000\,b^5c^2d^3 + 1102250152828125000000\,b\,c^2 d^3e^2 \\
& - 1785977460650880000000\,b^4c\,d^3f- 18582594433593750000\,b\,c^3 d\,e^3 + 553815035156250000000\,b^2c\,d^2e^3 - 132600864276720000000\,b^5 c^3d\,e  \\
& + 26161779208753125000000\,b\,c^2d^4f+ 282228536991744000000\,b^6 c\,d^2e - 497533299649804687500\,b^2 c^4d^2e - 21180712531289062500\,b^3c^4e\,f  \\
& - 23597441565143040000\,b^7c\,e\, f + 1309716804477312000000\,b^5{ d}^2e\,f - 70996261920117187500\,b^2 c^5d\,f - 20301909791625000000\,b^3 c^3d^2f \\
& + 47673706069335937500\,b\, c^6d\,e+ 1824711043695750000000\,b^3c^2d^3e - 14534321782640625000000\, b\,c\,d^5e + 43001246777343750000\,c^2d\,e^3f  \\
& + 182886113012695312500\,c^5d\,e\,f+ 12458939062500000000\,b\,c\,e^4f - 317507104115712000000\,b^5c\,d\,f^2 - 526406842089843750000\,b\,c^4d\,f^2  \\
&  - 134562534601420800000\,b^6 c^2 d\,f + 43655865760080000000\,b^5c\,e^2f - 8187922968750000000\,b\,c^4e^2f - 3014525999362500000000\,b^3d^2e^2 f  \\
& + 2113358174156250000000\,b\,d^2e^3f + 16929024323970937500000\, b^2c^2d^2f^2 + 2093420832890625000000\,c\,d^3e^2f + 280368861547851562500\,f^4e\,c\\
& - 8440672798215000000000\, b^2d\,e\,f^3 - 5262924058198242187500\,c^2d^2 e\,f^2 - 26161779208753125000000\,b\,d^3e\,f^2 + 23254914852225000000000\,b\, c\,d^2f^3  \\
& - 3749222721384000000\,b^6c^5 + 968815905029296875\,b^2c^8 + 9808823144531250000\,c^7{ d}^2 - 65051289843750000000\,f^2e^2c^3 \\
& + 1339843750000000000\,e^5 c^2 - 128892537789257812500\,b^2c^3e\, f^2 + 863536093567382812500\,b\, c^2e\,f^3 - 144501010939772928\,b^{10}c^2 \\
& \left.  - 1326295898437500000\,e\,c^8+ 10830575555419921875\,f^2c^6 - 56976855468750000000\,f^2e^4   \right)    \\
A_0 & =  3^4 \cdot 5^3   \,  \left(+733534659844875000000\,b^2c\,d^2e\,f^2  -48264980313867187500\,b\,c^3d\,e\, f^2- 3537852491250000000\,b^2 c\,d\,e^3f \right. \\
& - 245037899091037500000\,b^3c^2 d^2e\,f +53368632376980480000\,b^6c\,d\,e\,f+6304788593085937500\,b^2c^4d\,e\,f +2260894499521875000000\,b\,c\,d^4e\,f\\
& +12974633789062500\,c^{10} - 595795183593750000\,b\,c^8d -  1076935425804000000\,b^5c^5d+252876769144602624\,b^9c^2d  - 14469724796940000000\, b^4c\,d\,e^2f\\
& +367569090000000000\,b^2d^2e^4+9119299206445312500\,b^2c^6d^2+1675731938388480000\,b^6c^3d^2 - 54415270335431250000\,b^3c^4d^3+104182018537968000000\,b^4c^2d^4\\
& +93480771081600000\,b^6c^4e - 129317694433593750\,b^2c^7e - 1041820185379680000000\,b^3d^5e+1017211289062500000\,c^6d^2e + 51837087290625000000\,b\,c^2d^2e^2f \\
& +7586303074338078720\,b^8d^2e - 104702981241796875\,b^4c^4e^2 - 34870250374848000000\,b^6d^2e^2+39697461720000000000\,b^4d^2e^3 \\
& +125019882421875000\,b^2c^4e^3- 1369573279593750000\,b^3c^6f+ 2704690146170880000\,b^7c^3f - 124054567875000000000\,d^3e^3f - 27840669221998080000\,b^6c^2f^2 \\
& +2187911250000000000\,c^3d^2e^3+14224318264383897600\,b^7d\,f^2+6953147766386718750\,b^2c^5f^2+222692981458007812500\,c^4d^2f^2+2421516357421875000\,c^2e^3f^2 \\
& +3391341749282812500000\,d^4e\,f^2+1928463068847656250\,c^5e\,f^2 - 5767588008984375000\,b\,c^4f^3+129648734180582400000\,b^5c\,f^3 \\
& - 5087012623924218750000\,c\,d^3f^3 - 520910092689840000000\,b^3d^2f^3 + 167473666631250000000\,b^2e\,f^4+100530991698046875\,b^4c^7 \\
& - 99825278824857600\, b^8c^4+3488378906250000\,c^7e^2 - 16748046875000000\,c^4e^4+1134472500000000000\,d^2e^5 - 231515596751040000000\,b^4f^4               \\
& - 25233197539306640625\,c^3f^4+1695670874641406250000\,d\,f^5 - 30286759735107421875\,e^2f^4 - 4426806867259392000\,b^7c^2d\,e+4711911063398437500\,b^3c^5d\,e  \\
& - 33992369364240000000\,b^4c^3d^2e - 38735405887500000000\,b\,c^4d^3e+61737492466944000000\,b^5c\,d^3e+411267475627312500000\,b^2c^2d^4e            \\
& +154325882812500000\,b\,c^5d\,e^2+2198388719394000000\,b^5c^2d\,e^2 - 10588943472187500000\, b^2c^3d^2e^2+164347491520800000000\,b^3c\,d^3e^2       \\
& - 4782117619125000000\,b^3c^2d\,e^3 - 41351522625000000000\,b\,c\,d^3e^3+417323812500000000\,b\,c^2d\,e^4 - 3793151537169039360\,b^8c\,d\, f       \\
& +5780233979730000000\,b^4c^4d\, f+98903453635546875000\,b\,c^5d^2 f - 23813032808678400000\,b^5c^2d^2f - 842990663393156250000\,b^2c^3 d^3f        \\
& +1562730278069520000000\,b^3c\, d^4f+42383803710937500\,b\,c^6e\,f - 844705274742000000\,b^5c^3e\, f - 119624047593750000000\,c^3d^3e\, f           \\
& - 578788991877600000000\,b^4d^3e\, f+1553266568601562500\,b^3c^3e^2f - 11395371093750000\,b\,c^3e^3 f - 53591573322000000000\,b^2d^3e^2f           \\
& - 3486983554687500000\,c^4d\,e^2 f - 2127135937500000000\,c\,d\,e^4 f+8254059283968750000\,b^3c^3d\, f^2 - 3730475924211093750000\,b\,c^2d^3f^2     \\
& +659819450740464000000\,b^4c\, d^2f^2+3508972062750000000\,b^4c^2e\,f^2 - 191386226647526400000\,b^5d\,e\,f^2+154968966189450000000\,b^3d\,e^2f^2    \\
& - 8247228898535156250\,b^2c^2e^2f^2+209342083289062500000 c \,d^2e^2f^2 - 12060860765625000000\,b\,d\,e^3f^2 - 425263489195781250000\,b^2c^2d\,f^3  \\
& - 4845881557617187500\,b\,c\,e^2f^3+678268349856562500000\,b\,d^2e\,f^3 - 133330969713867187500\,c^2d\,e\, f^3 - 33494733326250000000\,b^3c\,e f^3    \\
& \left.  - 3051633867187500000\,c^7d\,f+791313074832656250000\,b\,c\,d\,f^4 \right)    \\
\end{split}
\end{equation*}

\end{tiny}

\end{document}